\documentclass[12pt]{article}
\usepackage[all]{xy}
\usepackage[dvips]{graphicx}
\usepackage{enumerate,hyperref}
\usepackage{ifthen,color}
\usepackage{amsmath,amsfonts,amssymb,amsthm,dsfont}
\usepackage{mathrsfs}
\usepackage[english]{babel}

\newtheorem{thm}{Theorem}[section]
\newtheorem{prop}[thm]{Proposition}
\newtheorem{coro}[thm]{Corollary}

\newtheorem{exa}[thm]{Example}
\newtheorem{lem}[thm]{Lemma}

\newtheorem{de}[thm]{Definition}

\newcommand{\comp}{ \,{\scriptstyle \stackrel{\circ}{}}\, }

\newcommand{\Acal}{\mathcal{A}}

\newcommand{\Bcal}{\mathcal{B}}

\newcommand{\Hcal}{\mathcal{H}}
\newcommand{\calH}{\mathcal{H}}

\newcommand{\RN}{ \mathbb{R} } 
 
\newcommand{\CN}{ \mathbb{C} } 
 
\newcommand{\HN}{\mathbb{H}}
\newcommand{\NN}{\mathbb{N}}

\newcommand{\eps}{\epsilon}

\newcommand{\mapmap}{\textbf{map}}

\newcommand{\MatA}{\text{M}_{\Acal}}

\newcommand{\EndA}{\mathcal{R}_{\Acal}}

\newcommand{\FunA}{\text{C}_{\Acal}(\Acal)}
\newcommand{\FunB}{\text{C}_{\Bcal}(\Bcal)}
\newcommand{\FunAj}{\text{C}_{\Acal_j}(\Acal_j)}

\newcommand{\MM}{\mathbf{M}}
\newcommand{\Acalx}{\mathcal{A}^{\times}}
\newcommand{\Acalzd}{\mathbf{zd}(\mathcal{A})}

\newcommand{\ds}{\displaystyle}

\newcommand{\Abound}{m_{\Acal}}

\newcommand{\zd}{\mathbf{zd}}
\newcommand{\Ld}{\mathrm{Ld}}

\begin{document}

\title{{\sc Introduction to the Theory of $\Acal$-ODEs}}
\author{Nathan BeDell\footnote{author was undergraduate at Liberty University during the creation of this work} \space\space\space\space\space\space\space James S. Cook \\ nbedell@tulane.edu \space\space jcook4@liberty.edu}
\maketitle
\begin{abstract}
We study the theory of ordinary differential equations over a commutative finite dimensional real associative unital algebra $\mathcal{A}$. We call such problems $\mathcal{A}$-ODEs. If a function is real differentiable and its differential is in the regular representation of $\mathcal{A}$ then we say the function is $\mathcal{A}$-differentiable. In this paper, we prove an existence and uniqueness theorem, derive Abel's formula for the Wronskian and establish the existence of a fundamental solution set for many $\mathcal{A}$-ODEs. We show the Wronskian of a fundamental solution set cannot be a divisor of zero. Three methods to solve nondegenerate constant coefficient $\mathcal{A}$-ODE are given. First, we show how zero-divisors complicate solution by factorization of operators. Second, isomorphisms to direct product are shown to produce interesting solutions. Third, our extension technique is shown to solve any nondegenerate $\mathcal{A}$-ODE; we find a fundamental solution set by selecting the component functions of the exponential on the characteristic extension algebra. The extension technique produces all of the elementary functions seen in the usual analysis by a bit of abstract algebra applied to the appropriate exponential function. On the other hand, we show how zero-divisors destroy both existence and uniqueness in degenerate $\mathcal{A}$-ODEs. We also study the Cauchy Euler problem for $\Acal$-Calculus and indicate how we may solve first order $\mathcal{A}$-ODEs.
\end{abstract}

\section{introduction and overview}
We use $\Acal$ to denote a real unital associative algebra of finite dimension. Elements of $\Acal$ are known as $\Acal$-numbers. We study calculus where real numbers have been replaced by $\Acal$-numbers. The resulting calculus we refer to as $\Acal$-calculus. Our typical goal is to find theorems which apply to as large a class of real commutative associative algebras as possible. In this paper we study the elementary theory of ordinary differential equations over $\Acal$. In particular, this means the differential equation, or system of differential equations, involve a set of dependent variables all of which depend on a single independent $\Acal$-variable. We call such differential equations $\Acal$-ODEs. \\

\noindent
In this paper we study: existence and uniqueness, first order problems, constant coefficient $n$-th order problems and Cauchy Euler problems over $\Acal$. In each topic we find either results or proofs require modification from the standard results over $\RN$ or $\CN$. Essentially, the existence of zero divisors forces nuances before unseen over a field and the submultiplicativity of the norm complicates the analysis. \\

\noindent
The main results of this paper include the existence and uniqueness Theorem \ref{thm:existenceofsystems} which provides the cornerstone for analysis of $\Acal$-ODEs. We also found many of the usual theorem for linear systems naturally generalize. For example, the Wronskian is still useful in that Abel's formula can be derived in the $\Acal$-calculus. However, instead of requiring the Wronskian be nonzero we must generalize to insist the Wronskian not be a divisor of zero in $\Acal$. That is, the Wronskian of linearly independent functions is a unit-valued in $\Acal$. This is simply the natural consequence of working in an $\Acal$-module as opposed to the standard theory where solutions form a vector space. We find seemingly strange results due to zero-divisors; it is possible to have distinct $\alpha_1, \alpha_2$ for which $e^{\alpha_1 x}$ and $e^{\alpha_2 x}$ are linearly dependent. This is impossible over a field. \\

\noindent
Section \ref{sec:firstorderAODEs} outlines how the usual elementary methods for solving first order differential equations abstract to the $\Acal$-calculus without much difficulty. However, when we examine the real PDE content of such problems we find nonlinear, coupled, systems of PDEs which are somehow solved by doing elementary calculus in $\Acal$. We hope this section helps the reader appreciate that the apparent simplicity of $\Acal$-calculus is a mask for something far less simple at the component level. The deeper question we would love to answer (but cannot at this time) is when we can solve the inverse problem; given a set of real PDEs what algebra $\Acal$ (if any) allows us to reformulate the system as an $\Acal$-ODE? \\

\noindent
We provide a calculational frame work to solve {\bf any} constant coefficient $\Acal$-ODE in Section \ref{sec:nondegencc}. While the general ideas have been known since the time of Euler, we think the method of Section \ref{sec:stupidsolutions} is new to the literature. In particular, we find it fascinating that the method provides a computational method which derives all the usual cases faced over $\RN$ in one sweeping algebraic method. Of course, the Laplace Transform also allows such a simplification, but, our method involves just a bit of abstract algebra and a natural chain rule. We find the special functions of the characteristic algebra {\bf always} provide a solution set in the nondegenerate case. This includes all the familiar elementary functions and a host of new functions. These new functions somehow spring into existence from the appearance of zero-divisors. For example, see Equation \ref{eqn:generalsol} where we see that roots which differ by a zero divisor are {\it almost} repeated. \\

\noindent
Degenerate $\Acal$-ODEs are not widely studied to our knowledge. In Section \ref{sec:degenerate} we were able to use algebra found in \cite{bedell1} to analyze zero-divisors interaction with linear operators over $\Acal$. We were pleased to provide two systematic families which illustrate how zero divisors destroy the usual theory. We found solution sets with infinite rank, and initial value problems for which no solution could be found. \\

\noindent
Finally, we are unaware of other works which treat the Cauchy-Euler problem in the generality we consider in Section \ref{sec:cauchyeuler}. We hope the reader is amused by the formula for the hyperbolic square root function given in Equation \ref{eqn:hyperbolicroot}. It seems likely this can be found in the literature, but, we are currently unaware of a reference.

\section{introduction to $\Acal$-calculus} \label{sec:reviewofAcalculus}

\noindent
In a nutshell, the study of $\Acal$-calculus is the study of calculus where real numbers have been replaced by numbers in an algebra $\Acal$. For example, complex analysis could be termed $\CN$-calculus, or the usual real calculus is $\RN$-calculus. Our focus is on the case $\Acal$ is associative, finite dimensional, and commutative. Some older works which align closely with our general methods are \cite{ward1940}, \cite{wagner1948} and in some sense \cite{ketchum}. We cannot hope to provide a complete history here, but, a better sketch is given in \cite{cookAcalculusI}.

\subsection{algebra and the regular representations}

\noindent
We say\footnote{we write $(\Acal, \star)$ to emphasize the pairing where helpful} $\Acal$ is an {\bf algebra} if $\Acal$ is a finite-dimensional real vector space paired with a function $\star: \Acal \times \Acal \rightarrow \Acal$ which is called {\bf multiplication}. In particular, the multiplication map satisfies the properties below: 
\begin{enumerate}[{\bf (i.)}]
\item {\bf bilinear:} $ (cx + y ) \star z = c( x \star z)+ y \star z$ and $x \star (cy+ z) = c(x \star y)+ x \star z $ for all $x,y,z \in \Acal$ and $c \in \RN$,
\item {\bf associative:} for which $x \star (y  \star z) = (x \star y)  \star z$  for all $x,y,z \in \Acal$ and,
\item {\bf unital:} there exists $\mathds{1} \in \Acal$ for which $\mathds{1} \star x = x$ and $x \star \mathds{1} = x$.
\end{enumerate}
We say $x \in \Acal$ is an {\bf $\Acal$-number}. If $x \star y = y \star x$ for all $x,y \in \Acal$ then $\Acal$ is {\bf commutative}. \\

\noindent
A linear transformation $T: \Acal \rightarrow \Acal$ is {\bf right $\Acal$ linear} if $T(x \star y) = T(x) \star y$ for all $x, y \in \Acal$. For example, $L_x(y) = x \star y$ defines a right-$\Acal$-linear map. We say the set $\EndA$ of all right $\Acal$ linear transformations forms the {\bf regular representation} of $\Acal$. Since $\Acal$ is unital the regular representation is isomorphic to $\Acal$. The isomorphism from $\Acal$ to $\EndA$ is given by $\mapmap(x) = L_x$ and we denote $\mapmap^{-1} = \#$ where $\#(T) = T(1)$. The idea here is that $\#(T)$ provides the $\Acal$ number which corresponds to $T$. If $\beta$ is a basis for $\Acal$ then the {\bf matrix regular representation} of $\Acal$ with respect to $\beta$ is
\begin{equation}
 \MatA(\beta) = \{ [T]_{\beta, \beta} \ | \ T \in \EndA \}. 
\end{equation} 
where $[T]_{\beta,\beta}$ denotes the matrix of $T$ with respect to the basis $\beta$. In the case $\Acal = \RN^n$ we may forgo the $\beta$ notation and write
\begin{equation}
 \MatA = \{ [T] \ | \ T \in \EndA \}
 \end{equation}
for the {\bf regular representation} of $\Acal$. There is a natural isomorphism of $\Acal$ and $\MatA$: If $\beta = \{ v_1 , \dots , v_n \}$ is a basis for $\Acal$ where $v_1 = \mathds{1}$ then 
\begin{equation}
\MM( x) = \left[ [x]_{\beta}| [x \star v_2]_{\beta}| \cdots | 
[x \star v_n]_{\beta} \right]
\end{equation}
where $[x]_{\beta}$ is the coordinate vector of $x$ with respect to $\beta$. In many applications we consider the case $\Acal = \RN^n$ with $\beta = \{ e_1, \dots , e_n \}$ the usual standard basis such that $e_1 = \mathds{1}$. Given these special choices we obtain much improved formula
\begin{equation}
 \MM(x) = [x| x \star e_2| \cdots | x \star e_n]. \end{equation}
We sometimes use {\it juxtaposition} in place of $\star$; $x \star y = xy$. For example:

\begin{exa} \label{Ex:number2}
The {\bf hyperbolic} numbers are given by $\calH = \RN \oplus j\RN$ where $j^2=1$. Identifying $e_1=1$ and $e_2=j$ we have $a+bj = [a,b]^T$. Moreover, 
\begin{equation}
(a+bj)e_2 = (a+bj)j = aj+b = [b,a]^T. 
\end{equation}
Therefore, $ \mathbf{M}(a+bj) = \left[ \begin{array}{cc} a  & b  \\  b & a \end{array} \right]$ is a typical matrix in $\text{M}_{\calH}$.
\end{exa}

\noindent
We say $x \in \Acal$ is a {\bf unit} if there exists $y \in \Acal$ for which $x \star y = y \star x = \mathds{1}$. The set of all units is known as the {\bf group of units} and we denote this by $\Acalx$. We say $a \in \Acal$ is a {\bf zero-divisor} if $a \neq 0$ and there exists $b \neq 0$ for which $a \star b = 0$ or $b \star a =0$. Let $\Acalzd = \{ x \in \Acal \ | \ x=0 \ \text{or $x$ is a zero-divisor} \}$.

\begin{exa} \label{prop:isomorphismhyperbolic}
 and $\mathbf{zd}( \calH ) = \{ a+bj \ | \ a^2=b^2 \}$ whereas $\calH^{\times} = \{ a+bj \ | \ a^2 \neq b^2 \}$. The reciprocal of an element in $\calH^{ \times}$ is simply
\begin{equation} 
\frac{1}{a+bj} = \frac{a-bj}{a^2-b^2}. 
\end{equation}
The above follows from the identity $(a+bj)(a-bj) = a^2-b^2$ given $a^2-b^2 \neq 0$. Let $\Bcal = \RN \times \RN$ with $(a,b)(c,d) = (ac, bd)$ for all $(a,b),(c,d) \in \Bcal$. We can show that
\begin{equation} \Psi(a,b) = a\left(\frac{1+j}{2}\right)+b\left(\frac{1-j}{2}\right) \qquad \& \qquad
\Psi^{-1}(x+jy) = (x+y,x-y) 
\end{equation} 
provide an isomorphism of $\calH$ and $\RN \times \RN$. In \cite{cookAcalculusI} an examples are given which show how this isomorphism can be used to solve the quadratic equation in $\Hcal$ and to derive  d'Alembert's solution to the wave equation.
\end{exa}

\subsection{submultiplicative norms}
The division algebras $\RN, \CN$ and $\HN$ can be given a {\bf multiplicative norm} where $\| x \star y \| = \| x \| \, \| y \|$. Generally we can only find {\bf submultiplicative norm}. 

\begin{exa}
 If $\Hcal$ is given norm $\| x+jy \| = \sqrt{x^2+y^2}$ then $\| zw \| \leq \sqrt{2} \|z \| \, \|w \|$.
\end{exa} 

\noindent
If $\Acal$ is an algebra over $\RN$ with basis $\{ v_1, \dots , v_n \}$ then define {\bf structure constants} $C_{ijk}$ by $v_i \star v_j = \sum_{k=1}^n C_{ijk} v_k$ for all $1 \leq i,j  \leq n$. For proof of what follows see \cite{cookAcalculusI}.   

\begin{thm} \label{thm:submultiplicative}
\textbf{(submultiplicative norm)} If $\Acal$ is an associative $n$-dimensional algebra over $\RN$ then there exists a norm $|| \cdot ||$ for $\Acal$ and $\Abound >0$ for which $|| x \star y|| \leq \Abound ||x|| ||y||$ for all $x,y \in \Acal$. Moreover, for this norm we find $\Abound = \mathbf{C}(n^2-n+1) \sqrt{n}$ where $\mathbf{C} = \text{max} \{ C_{ijk} \ | \ 1 \leq i,j,k \leq n \}$.
\end{thm}

\begin{coro} \label{prop:ineqnpower}
If $\| x \star y \| \leq \Abound \| x \| \|y \|$ for $x,y \in \Acal$ then $\| z^n \| \leq \Abound^n \| z \|^n$ for each $z \in \Acal$ and $n \in \NN$.
\end{coro}

\begin{coro} \label{thm:quotientinequality}
Suppose $\Abound >0$ is a real constant such that $|| x \star y|| \leq \Abound ||x|| ||y||$ for all $x,y \in \Acal$. If $b \in \Acalx$ and $a \in \Acal$ then $ \frac{||a||}{||b||} \leq \Abound \, \big{|}\big{|} \frac{a}{b} \big{|}\big{|}$. 
\end{coro}

\subsection{differential calculus on $\Acal$}
The definition of differentiability with respect to an algebra variable is open to some debate. There seem to be two main approaches:
\begin{quote}
\begin{enumerate}[{\bf D1:}]
\item define differentiability in terms of an algebraic condition on the differential,
\item define differentiability in terms of a deleted-difference quotient.
\end{enumerate}
\end{quote}
In \cite{cookAcalculusI} it is shown that these definitions are interchangeable on an open set in the context of a commutative semisimple algebra. However, it is also shown that in there exist {\bf D1} differentiable functions which are nowhere {\bf D2}. Hence, we prefer to use {\bf D1} as it is more general. Following \cite{cookAcalculusI} we define differentiability with respect to an algebra variable as follows:

\begin{de} \label{defn:Adiff}
Let $U \subseteq \Acal$ be an open set containing $p$. If $f: U \rightarrow \Acal$ is a function then we say $f$ is {\bf $\Acal$-differentiable at $p$} if there exists a linear function $d_pf \in \EndA$ such that
\begin{equation} \label{eqn:frechetquotAcal}
 \lim_{h \rightarrow 0}\frac{f(p+h)-f(p)-d_pf(h)}{||h||} = 0. 
\end{equation}
 \end{de}

\noindent
In other words, $f$ is $\Acal$-differentiable at a point if its differential at the point is a right-$\Acal$-linear map. Equivalently, given a choice of basis, $f$ is $\Acal$-differentiable if its Jacobian matrix is found in the matrix regular representation of $\Acal$. If $\Acal$ has basis $\beta = \{ v_1, \dots , v_n\}$ has coordinates $x_1, \dots , x_n$ then $d_pf(e_j) = \frac{\partial f}{\partial x_j}(p)$. Suppose $v_1=\mathds{1}$ then $ d_pf(1) = \frac{\partial f}{\partial x_1}(p)$. Observe right linearity of the differential indicates $d_pf(v_j) = d_pf(\mathds{1} \star v_j) = d_pf(1) \star v_j$ hence for each $p$ at which $f$ is $\Acal$ differentiable we find:
\begin{equation}
\frac{\partial f}{\partial x_j} (p)= \frac{\partial f}{\partial x_1}(p) \star v_j.
\end{equation}
These are the {\bf $\Acal$-Cauchy Riemann Equations}. There are $n-1$ equations in $\Acal$ which amount to $n^2-n$ scalar equations. If the $\Acal$-CR equations hold for a continuously differentiable $f$ at $p$ then we have that $d_pf \in \EndA$.  \\

\noindent
Next we wish to explain how to construct the derivative function $f'$ on $\Acal$.  We are free to use the isomorphism between the right $\Acal$ linear maps and $\Acal$ as to define the {\it derivative at a point} for via $f'(p) = \#(d_pf)$. This is special to our context. In the larger study of real differentiable functions on an $n$-dimensional space no such isomorphism exists and it is not possible to identify arbitrary linear maps with points.

\begin{de} \label{defn:derivative}
Let $U\subseteq \Acal$ be an open set and $f: U  \rightarrow \Acal$ an $\Acal$-differentiable function on $U$ then we define $f': U \rightarrow \Acal$ by $f'(p) = \# (d_pf) $ for each $p \in U$.
\end{de}

\noindent
Equivalently, we could write $f'(p) = d_pf( \mathds{1})$ since $\#(T) = T( \mathds{1})$ for each $T \in \EndA$.
Many properties of the usual calculus hold for $\Acal$-differentiable functions.  

\begin{prop}
For $f$ and $g$ both $\Acal$-differentiable at $p$,
\begin{quote}
 \begin{enumerate}[{\bf (i.)}]
\item $ \ds (f+g)'(p) = f'(p)+g'(p)$,
\item for constant $c \in \Acal$, $ \ds (c \star f)'(p) = c \star f'(p)$,
\item given $\Acal$ is commutative, $ (f \star g)'(p) = f'(p) \star g(p)+ f(p) \star g'(p)$,
\item $ (f \comp g)'(p) = f'(g(p)) \star g'(p), $
\item if $f(\zeta) = \zeta^n$ for some $n \in \NN$ then $f'(\zeta) = n \zeta^{n-1}$. 
\end{enumerate}
\end{quote}
\end{prop}

\noindent
If $\Acal$ is not commutative then the product of $\Acal$-differentiable functions need not be $\Acal$-differentiable. In \cite{cookAcalculusI} an example is given where $f,g$ and $f \star g$ are $\Acal$-differentiable yet $g \star f$ is not $\Acal$-differentiable. \\

\noindent
We are also able to find an $\Acal$-generalization of Wirtinger's calculus. In \cite{cookAcalculusI} we introduce conjugate variables $\bar{\zeta}_2, \dots, \bar{\zeta}_n$ for $\Acal$ and find for commutative algebras if $f: \Acal \rightarrow \Acal$ is $\Acal$-differentiable at $p$ then $\ds \frac{\partial f}{\partial \overline{\zeta}_j} = 0$ for $j=2, \dots , n$. In other words, another way we can look at $\Acal$-differentiable functions is that they are functions of $\zeta$ alone.  \\  

\noindent
The theory of higher derivatives is also developed in \cite{cookAcalculusI}.

\begin{de} \label{defn:higherderivative}
Suppose $f$ is a function on $\Acal$ for which the derivative function $f'$ is $\Acal$-differentiable at $p$ then we define
$ f''(p) = (f')'(p)$. Furthermore, supposing the derivatives exist, we define $f^{(k)}(p) = (f^{(k-1)})'(p)$ for $k =2,3, \dots$. 
\end{de}

\noindent
Naturally we define functions $f'', f''', \dots, f^{(k)}$ in the natural pointwise fashion for as many points as the derivatives exist. Furthermore, with respect to $\beta = \{ v_1, \dots , v_n \}$ where $v_1 = \mathds{1}$, we have $f'(p) = d_pf(\mathds{1}) = \frac{\partial f}{\partial x_1}(p)$.  Thus, $f' = \frac{\partial f}{\partial x_1}$. Suppose $f''(p)$ exists. Note,
\begin{equation}
 f''(p) = (f')'(p) = \#( d_p f'(\mathds{1}) )= \frac{\partial f'}{\partial x_1}(p) = \frac{\partial^2f}{\partial x_1^2}(p). 
\end{equation}
Thus, $f'' = \frac{\partial^2 f}{\partial x_1^2}$. By induction, we find if $f^{(k)}$ exists then $f^{(k)} = \frac{\partial^k f}{\partial x_1^k}$. Furthermore, if $f: \Acal \rightarrow \Acal$ is $k$-times $\Acal$-differentiable then 
\begin{equation}
\frac{\partial^k f}{\partial x_{i_1}\partial x_{i_2} \cdots \partial x_{i_k}} = \frac{\partial^k f}{\partial x_1^k} \star v_{i_{1}} \star v_{i_{2}} \star \cdots \star v_{i_{k}}.
\end{equation}

\noindent
The Theorem below gives us license to convert equations in $\Acal$ to partial differential equations which every component of an $\Acal$-differentiable function must solve!

\begin{thm}
Let $U$ be open in $\Acal$ and suppose $f: U \rightarrow \Acal$ is $k$-times $\Acal$-differentiable. If there exist $B_{i_1i_2\dots i_k} \in \RN$ for which $\sum_{i_1i_2\dots i_k} B_{i_1i_2\dots i_k}v_{i_1} \star v_{i_2} \star \cdots \star v_{i_k} = 0$ then 
\begin{equation}
\sum_{i_1i_2\dots i_k} B_{i_1i_2\dots i_k}\frac{\partial^k f}{\partial x_{i_1} \partial x_{i_2} \cdots  \partial x_{i_k} }=0. 
\end{equation}
\end{thm}

\begin{exa}
Since $i^2=-1$ in $\CN$ it follows for $z=x+iy$ that complex differentiable $f$ have $f_{yy}=-f_{xx}$. Setting $f=u+iv$ we find $u_{xx}+u_{yy}=0$ and $v_{xx}+v_{yy}=0$. In other words, complex differentiable functions solve Laplace's equation $u_{xx}+u_{yy}=0$ because $i^2+1=0$. Likewise, in $\Hcal = \RN \oplus j \RN$ we find solutions to the one-dimensional wave equation since $j^2-1=0$ implies $u_{xx}-u_{yy}=0$ and $v_{xx}-v_{yy}=0$ for $u+jv$ a $\Hcal$-differentiable function. 
\end{exa}

\noindent
One may ask when a given set of real PDEs appears as the $\Acal$-CR equations or one of their differential consequences. For example, Ward showed in \cite{ward1952} that if we are given an appropriate set of real PDEs then we can find $\Acal$ for which those PDEs are the $\Acal$-CR equations. The general problem of ascertaining if a given set of PDEs is consistent with the function theory for a given algebra $\Acal$ is in our estimation a difficult and open question. 

\subsection{integral calculus on $\Acal$}

\noindent Integration along curves in $\Acal$ is defined in \cite{cookAcalculusI} in much the same fashion as $\CN$. If $\zeta: [t_o,t_1] \rightarrow \Acal$ is differentiable parametrization of a curve $C$ and $f$ is continuous near $C$ then 
\begin{equation}
 \ds \int_C f ( \zeta) \star d\zeta = \int_{t_o}^{t_f} f( \zeta (t)) \star \frac{d\zeta}{dt} \, dt. \end{equation}

\begin{thm} \label{thm:subML}
Let $C$ be a rectifiable curve with arclength $L$. Suppose $||f(\zeta) || \leq M$ for each $\zeta \in C$ and suppose $f$ is continuous near $C$. Then 
$$ \bigg{|}\bigg{|} \int_C f( \zeta) \star d\zeta \bigg{|}\bigg{|} \leq \Abound ML $$
where $\Abound$ is a constant such that $|| z \star w || \leq \Abound ||z|| \, ||w||$ for all $z,w \in \Acal$.
\end{thm}
\noindent
First and Second Fundamental Theorems of Calculus, discussion of exact and closed forms, basic topological theorems on line-integrals, Cauchy's Integral Theorem are all found in $\cite{cookAcalculusI}$.

\subsection{further background}
\noindent
This paper follows from a number of papers which at least one of the authors has participated. In particular, the trouble of submultiplicative norms is more forcefully seen in \cite{cookfreese} where the theory of convergence and divergence for power series in $\Acal$ is studied in depth. For this paper, we simply need the existence of series expansions and the properties of elementary functions. Further details about definitions of sine, cosine and the exponential as well as the properties for an arbitrary commutative, associative finite dimensional real algebra see \cite{cookfreese}, \cite{bedell3} and \cite{bedell1}. Properties of logarithms over many algebras are studied in \cite{bedell2}. The reader may find the many explicit examples in \cite{cook} a helpful supplement to our current work.

\section{existence and uniqueness for $\Acal$-ODEs} \label{sec:existenceanduniqueness}
We assume $\Acal$ is commutative throughout this Section. 
If $\vec{g}: \Acal^m \rightarrow \Acal^k$ then $\vec{g} = (g_1, \dots, g_k)$ is $\Acal$-differentiable if each component function $g_i$ is $\Acal$-differentiable in the sense that $\vec{g}$ is real differentiable and has a right-$\Acal$-linear differential. Let $\Acal^m$ have algebra variables $z_1,z_2, \dots , z_m$. If $\vec{g}$ is $\Acal$-differentiable at $p$ then define
\begin{equation}
\frac{\partial \vec{g}}{\partial z_i}(p) = d_p\vec{g}(e_i)
\end{equation}
where $e_1 = (1,0,\dots , 0), \dots , e_m = (0,\dots , 0, 1)$. Hence generally,
\begin{equation} \label{eqn:partialsetup}
d_p\vec{g}(h_1,h_2,\dots , h_m) = h_1\frac{\partial \vec{g}}{\partial z_1}(p)+
h_2\frac{\partial \vec{g}}{\partial z_2}(p)+
+ \cdots +
h_m\frac{\partial \vec{g}}{\partial z_m}(p).
\end{equation}
These partial derivatives are used in what follows. 

\begin{thm} \label{thm:existenceofsystems}
Let $I$ be compact and star-shaped in $\Acal$ and let $R = I \times \Acal^k$. Suppose $\vec{f}: R \rightarrow \Acal^k$ is $\Acal$-differentiable with $\big{|}\big{|} \frac{\partial  \vec{f}}{\partial y_i}(z, \vec{y}) \big{|}\big{|} \leq L$ for each $(z, \vec{y}) \in R$. Let $(z_o, \vec{w}_o) \in R$. The initial value problem $\frac{d\vec{y}}{dz} = \vec{f}(z, \vec{y})$ with $\vec{y}(z_o) = \vec{w}_o$ has a unique solution on $I$.
\end{thm}

\noindent
{\bf Proof:} we intend to define the solution as the limit function of a Picard iteration. Begin by setting $\vec{y}_o = \vec{w}_o$ and for $n=0,1,\dots $
\begin{equation}
\vec{y}_{n+1}(z) = \vec{w}_o+ \int_{[z_o,z]} \vec{f}(  \zeta, \vec{y}_n( \zeta)) \star d \zeta. 
\end{equation}
Since $I$ is star-shaped we know $[z_o ,z]  = \{ z_o+t(z-z_o) \ | \ 0 \leq t \leq 1 \} \subseteq I$ and consequently the integral is well-defined. In particular, we define for $\vec{f}= (f_1,f_2,\dots , f_k)$,
\begin{equation}
  \int_{[z_o,z]} \vec{f} \star d \zeta = 
\left( \int_{[z_o,z]} f_1 \star d \zeta, \int_{[z_o,z]} f_2 \star d \zeta, \dots , \int_{[z_o,z]} f_k \star d \zeta \right).
\end{equation}
Notice,
\begin{equation}\label{eqn:telescopesareneat}
 \vec{y}_o+ \sum_{j=0}^{n-1} \left[ \vec{y}_{j+1}-\vec{y}_j\right] = \vec{y}_o + \left[ \vec{y}_{1}-\vec{y}_o\right]+ \left[ \vec{y}_2-\vec{y}_1\right]+\cdots + \left[ \vec{y}_{n}-\vec{y}_{n-1}\right] =  \vec{y}_{n}. \end{equation}
Thus uniform convergence of $ \sum_{j=0}^{\infty} \left[ \vec{y}_{j+1}-\vec{y}_j\right]$ provides uniform convergence of $\{ \vec{y}_n\}$. We will show that $ \sum_{j=0}^{\infty} \left[ \vec{y}_{j+1}-\vec{y}_j\right]$ can be majorized over $I$ by a convergent series. Then, using  \cite{cookfreese}, we deduce the convergence of the series is uniform. \\

\noindent
Observe $I$ compact implies there exists $M>0$ for which $||f(\zeta, \vec{w}_o) || \leq M$ for all $ \zeta \in I$. Thus for $[z_o, z] \subset I$,
\begin{equation} \label{eqn:firstbound}
|| \vec{y}_1(z)- \vec{y}_o(z) || = \bigg{|}\bigg{|} \int_{[z_o,z]}\vec{f}(\zeta, \vec{w}_o) \star d\zeta \bigg{|}\bigg{|} \leq M\Abound||z-z_o|| 
\end{equation}
as the length of $[z_o,z]$ is simply $||z-z_o||$. We need a generalization of the mean value theorem for our current context to make further progress in the proof: 

\begin{lem} \label{thm:lemmaMVT}
With $\vec{f}$ and $R$ as in preceding discussion, there exists $l>0$ for which \\ $|| \vec{f}(\zeta, \vec{v})-\vec{f}(\zeta,\vec{w})|| \leq l|| \vec{v}-\vec{w}||$ for $(\zeta, \vec{v}),(\zeta, \vec{w}) \in R$.
\end{lem}

\noindent
{\bf Proof:} Notice, if $D\vec{f}$ denotes the Frechet derivative for $\vec{f}: R \rightarrow \Acal^k$ then Theorem 1 on page 73 of \cite{zorich} gives that 
\begin{equation} \label{eqn:Zorich}
 || \vec{f}(A+H) - \vec{f}(A)|| \leq ||H|| \, \text{sup}_{x \in [A,A+H]} \{ ||D\vec{f}(x)|| \} 
 \end{equation}
where we suppose $R$ is given norm by $||(z, \vec{y})|| = \sqrt{ ||z||^2+ ||y_1||^2+ \cdots + ||y_k||^2 }$ for each $(z,\vec{y}) \in R$ and $||D\vec{f}||$ denotes the operator norm defined by
\begin{equation}
 ||D\vec{f}(x)|| = \text{sup}_{||W||=1}(||D\vec{f}(x)(W)||)
\end{equation}
where $D\vec{f}(x)(W) =(y_1, \dots , y_k) \in \Acal^k$ has norm $||(y_1, \dots , y_k)|| = \sqrt{ ||y_1||^2+\cdots + ||y_k||^2}$
. Notice, from our initial discussion leading to Equation \ref{eqn:partialsetup},
\begin{equation} \label{eqn:differentialexpaspartials}
 (D\vec{f})(P)(w_o,\vec{w}) =  \frac{\partial \vec{f}}{\partial z}(P)w_o+ \frac{\partial \vec{f}}{\partial y_1}(P)w_1+\frac{\partial \vec{f}}{\partial y_2}(P)w_2+ \cdots + \frac{\partial \vec{f}}{\partial y_k}(P)w_k. 
 \end{equation}
If $P \in [A,A+H] \subset R$ where $H = (0,\vec{h})$ then $w_o=0$ whereas $w_i=h_i$ for $i=1,2,\dots, k$ so
\begin{equation}
 (D\vec{f})(P)(0,\vec{h}) =  \frac{\partial \vec{f}}{\partial y_1}(P)h_1+\frac{\partial \vec{f}}{\partial y_2}(P)h_2+ \cdots + \frac{\partial \vec{f}}{\partial y_k}(P)h_k. 
\end{equation}
We assumed $\big{|}\big{|} \frac{\partial  \vec{f}}{\partial y_i}(z, \vec{y}) \big{|}\big{|} \leq L$ for each $(z, \vec{y}) \in R$ hence by the triangle inequality and submultiplicativity of the norm on $\Acal$,
\begin{equation}
 ||(D\vec{f})(P)(0,\vec{h})|| \leq   \Abound ||h_1||L+\cdots +  \Abound ||h_k||L \leq kL \Abound ||\vec{h}||
\end{equation}
and $||D\vec{f}(P)|| \leq kL\Abound$ for $P \in [A,A+H] \subset R$ where $H = (0,\vec{h})$. Thus from \ref{eqn:Zorich} we find
\begin{equation} \label{eqn:ZorichII}
 || \vec{f}(A+H) - \vec{f}(A)|| \leq ||H||  kL \Abound 
 \end{equation}
 for $H = (0, \vec{h})$ with $[A,A+H] \subset R$. If $(\zeta, \vec{v}), (\zeta, \vec{w}) \in R$ then set $A = (\zeta, \vec{w})$ and $A+H = (\zeta, \vec{v})$ hence $H = (0, \vec{v}-\vec{w}) \in R$ and for $l = kL\Abound$ we find
 $|| \vec{f}(\zeta, \vec{v})-\vec{f}(\zeta,\vec{w})|| \leq l|| \vec{v}-\vec{w}||$. $\Box$ \\

\noindent
We now continue the proof of Theorem \ref{thm:existenceofsystems}. Let $l = kL\Abound$ and inductively suppose 
\begin{equation} \label{eqn:indhypothesis}
|| \vec{y}_n(z)- \vec{y}_{n-1}(z) ||   \leq \frac{Ml^{n-1}\Abound^n||z-z_o||^n}{ n!}
\end{equation}
for $[z_o,z] \subset I$. Notice Equation \ref{eqn:firstbound} gives the induction claim for $n=1$.  Consider, for $z \in I$,   
\begin{align} 
|| \vec{y}_{n+1}(z)- \vec{y}_{n}(z) || &= \bigg{|}\bigg{|} \int_{[z_o,z]}\left( \vec{f}(\zeta, \vec{y}_n(\zeta)) - \vec{f}(\zeta, \vec{y}_{n-1}(\zeta)) \right) \star d\zeta \bigg{|}\bigg{|} & \\ \notag
&\leq  \Abound l \int_{[z_o,z]} ||\vec{y}_n(\zeta)-\vec{y}_{n-1}(\zeta)|| \, ||  d\zeta ||  \ \ \ \text{( by Lemma \ref{thm:lemmaMVT}, )}\\ \notag
&\leq  \Abound l \int_{[z_o,z]} \frac{Ml^{n-1}\Abound^n||\zeta-z_o||^n}{n!} \, ||  d\zeta ||  \ \ \ \text{( by induction claim of \ref{eqn:indhypothesis}, )}\\ \notag
& = \frac{Ml^{n}\Abound^{n+1}}{n!}\int_{[z_o,z]} s^n \, ds  \\ \notag
&=  \frac{Ml^{n}\Abound^{n+1}}{n!}\cdot \frac{||z-z_o||^{n+1}}{n+1}
\end{align}
hence $||\vec{y}_{n+1}(z)- \vec{y}_{n}(z)|| \leq \frac{Ml^{n}\Abound^{n+1}||z-z_o||^{n+1}}{(n+1)!}$ and we find estimate \ref{eqn:indhypothesis} is true for all $n \in \NN$ by induction. Furthermore, if $s$ denotes the distance from $z_o$ to $z$ and $\beta = \Abound s l$ then we may reformulate the bound of \ref{eqn:indhypothesis} as
\begin{equation}
 || \vec{y}_n(z)- \vec{y}_{n-1}(z) ||   \leq \frac{M}{l} \cdot \frac{\beta^n}{n!}.
\end{equation}
Since $I$ compact we know there exists $s_o>0$ for which the distance $s=||z-z_o|| \leq s_o$. Let $\beta_o = \Abound l s_o$ and note that
\begin{equation}
 || \vec{y}_n(z)- \vec{y}_{n-1}(z) ||   \leq \frac{M}{l} \cdot \frac{\beta_o^n}{n!}
\end{equation}
for all $z \in I$. Since $\sum_{n=0}^{\infty}\frac{\beta_o^n}{n!} = e^{\beta_o}$ we have majorized the series $\sum_{n=0}^{\infty} || \vec{y}_n(z)- \vec{y}_{n-1}(z) ||$ on $I$. Thus, $\sum_{n=0}^{\infty} || \vec{y}_n(z)- \vec{y}_{n-1}(z) ||$ is uniformly convergent on $I$ and we deduce from Equation \ref{eqn:telescopesareneat} that $\{ \vec{y}_n \}$ converges uniformly to $\vec{y}_*$ on $I$. \\

\noindent
Let us examine why $\vec{y}_*$ is a solution to the initial value problem. First, note $\vec{y}_n(z_o)= \vec{w}_o$ and as uniform convergence implies pointwise convergence we have
\begin{equation}
 \vec{y}_*(z_o) = \left( \lim_{n \rightarrow \infty} \vec{y}_n \right)(z_o) =\lim_{n \rightarrow \infty} \left( \vec{y}_n (z_o) \right)  = \lim_{n \rightarrow \infty} \vec{w}_o = \vec{w}_o. 
\end{equation}
Second, to see $\vec{y}_*$ is a solution for $z \in I$, consider
\begin{equation} \vec{y}_*(z) = \lim_{n \rightarrow \infty} \vec{y}_n(z) = \lim_{n \rightarrow \infty} \left( \vec{w}_o+
 \int_{[z_o,z]} \vec{f}(  \zeta, \vec{y}_{n-1}( \zeta)) \star d \zeta \right).
\end{equation}
However, uniform convergence of $\{\vec{y}_n \}$ and continuity of $\vec{f}$ imply uniform convergence of $\{\vec{f}(  \zeta, \vec{y}_{n-1}( \zeta)) \}$ therefore we can exchange the order of integration and the limit to deduce
\begin{equation}
 \vec{y}_*(z) = \vec{w}_o+ \int_{[z_o,z]}  \left( \lim_{n \rightarrow \infty}\vec{f}(  \zeta, \vec{y}_{n-1}( \zeta)) \right) \star d \zeta = 
\vec{w}_o + \int_{[z_o,z]} \vec{f}(  \zeta, \vec{y}_*( \zeta))  \star d \zeta.
\end{equation}
Thus, $\frac{d\vec{y}_*}{dz} = \vec{f}(  z, \vec{y}_*( z))$ for each $z \in I$. \\

\noindent
Finally, to see the solution is unique, suppose $\vec{y}_{**}$ is a solution on $I$ of $\frac{d\vec{y}}{dz} = \vec{f}(z,\vec{y})$ with $\vec{y}_{**}(z_o) = \vec{w}_o$. Let $z \in I$, by Lemma \ref{thm:lemmaMVT}, $||\vec{f}(\zeta, \vec{y}_{**}(\zeta))- \vec{f}(\zeta, \vec{y}_{*}(\zeta)) || \leq l ||\vec{y}_{**}(\zeta)-\vec{y}_{*}(\zeta)||$ where $l>0$. Moreover,  $\Upsilon =\text{sup}\{ ||\vec{y}_{**}(\zeta)-\vec{y}_{*}(\zeta)|| \ | \ \zeta \in [z_o,z] \}$ provides a bound for $||\vec{y}_{**}(\zeta)-\vec{y}_{*}(\zeta)||$ on $[z_o,z]$ hence
\begin{equation} \label{eqn:stageone}
||\vec{y}_{**}(z)-\vec{y}_{*}(z)|| = 
\bigg{|}\bigg{|} \int_{[z_o,z]}\left( \vec{f}(\zeta, \vec{y}_{**}(\zeta))- \vec{f}(\zeta, \vec{y}_{*}(\zeta))\right) \star d\zeta \bigg{|}\bigg{|} \\ \leq l \cdot \Upsilon \cdot ||z-z_o||.
\end{equation}
Thus, by Lemma \ref{thm:lemmaMVT} and the estimate above,
\begin{equation}
||\vec{f}(\zeta, \vec{y}_{**}(\zeta))- \vec{f}(\zeta, \vec{y}_{*}(\zeta)) || \leq l ||\vec{y}_{**}(\zeta)-\vec{y}_{*}(\zeta)|| \leq l^2 \Upsilon ||\zeta-z_o||. 
\end{equation}
Thus,
\begin{equation} \label{eqn:stagetwo}
||\vec{y}_{**}(z)-\vec{y}_{*}(z)|| \leq l^2 \Upsilon \int_{[z_o,z]}||\zeta-z_o|| \, d\zeta = l^2 \Upsilon \int_{[z_o,z]} s \, ds  =\frac{l^1\Upsilon||z-z_o||}{2}.
\end{equation}
Continuing in the above fashion we find
\begin{equation}
 ||\vec{y}_{**}(z)-\vec{y}_{*}(z)|| \leq \frac{l^{n-1}\Upsilon||z-z_o||^n}{n!}. 
\end{equation}
for $n \in \NN$. As $n \rightarrow \infty$ we find $||\vec{y}_{**}(z)-\vec{y}_{*}(z)|| \rightarrow 0$ for each $z \in I$. Thus $\vec{y}_{**}(z)=\vec{y}_{*}(z)$ for each $z \in I$ and the proof of Theorem  \ref{thm:existenceofsystems} is complete. $\Box$ \\

\noindent
With the Theorem above in hand the remaining theory of linear $\Acal$-ODEs follows easily.

\begin{thm} \label{thm:existenceIVP}
Let $L = D^k+a_{k-1}D^{k-1}+\cdots + a_2D^2+a_1D+a_o$ where $a_o, a_1, \dots , a_{k-1}$ are $\Acal$-differentiable functions on a compact and star-shaped domain $I$ and $D = d/dz$. Also, suppose $g$ is an $\Acal$-differentiable function on $I$. The $k$-th order $\Acal$-ODE $L[y] = g$ with initial conditions $y(z_o)=y_o, y'(z_o)=y_1, \dots , y^{(k-1)}(z_o) = y_{k-1}$ for $z_o \in I$ has unique solution on $I$.
\end{thm}

\noindent
{\bf Proof:} the proof is by the usual reduction of order. Let $w_1=y, w_2=y', \dots , w_{k} = y^{(k-1)}$. Observe, $\frac{dw_j}{dz} = \frac{dy^{(j-1)}}{dz} = y^{(j)} = w_{j+1}$ for $j=1,2,\dots , k-1$. Observe $L[y]=0$ provides:
\begin{align}
 y^{(k)} &= g-a_{k-1}y^{(k-1)} - \cdots - a_2y''-a_1y'-a_oy  \\ \notag
 &= g-a_ow_1  -a_1w_2 - a_2w_3 - \cdots -a_{k-1}w_k. 
 \end{align}
Thus, as $w_{k}' = y^{(k)}$, we calculate the reduced system has a coefficient matrix which is a complementary matrix\footnote{or the transpose of a complementary matrix if you prefer} to the characteristic polynomial of the given $k$-th order $\Acal$-ODE,
\begin{equation}
 \frac{d\vec{w}}{dz} = A\vec{w} +\vec{b} \ \ \text{where} \ \ A = \left[ \begin{array}{ccccc} 
0 & 1  & \cdots & 0 & 0 \\
0 & 0  & \cdots & 0 & 0 \\
\vdots  & \vdots & \ddots & \vdots & \vdots \\
0 &  0 & \cdots & 0 & 1 \\
-a_o & -a_1  & \cdots & -a_{k-2} & -a_{k-1}
\end{array}\right] \ \ \& \ \ \vec{b} = \left[ \begin{array}{c} 0 \\ 0 \\ \vdots \\ 0 \\ g \end{array}\right]. 
\end{equation}
Notice $\vec{f}(z,\vec{w}) = A\vec{w}+\vec{b}$ is $\Acal$-differentiable since we suppose the coefficient functions $a_o, \dots , a_{k-1}$ and forcing term $g$ are $\Acal$-differentiable on $I$. Thus, by Theorem \ref{thm:existenceofsystems} we find a unique solution to $\frac{d\vec{w}}{dz} = A\vec{w} +\vec{b}$ for a given an initial condition vector $\vec{w}(z_o) =(y_o,y_1,\dots, y_k) \in \Acal^k$. By construction, $w_1=y$ of the solution provides the solution to the initial value problem $L[y]=g$ where $y(z_o)=y_o, y'(z_o)=y_1, \dots , y^{(k-1)}(z_o) = y_{k-1}$. $\Box$ \\

\noindent
The set of $\Acal$-differentiable functions has a natural $\Acal$-module structure. Hence define:

\begin{de} 
Let $I$ be a connected subset of $\Acal$. Suppose $f_j: I \rightarrow \Acal$ are functions. We say the set of functions $\{f_1,f_2,f_3, \dots , f_m\}$ are {\bf linearly independent (LI)} on $I$ if and only if for $c_1, \dots , c_m \in \Acal$
\[ c_1f_1(z)+c_2f_2(z)+ \cdots +c_mf_m(z) = 0 \]
for all $z \in I$ implies $c_1=c_2= \cdots =c_m=0$. Conversely, if $\{f_1,f_2,f_3, \dots , f_m\}$ are not linearly independent on $I$ then they are said to be {\bf linearly dependent} on $I$.  
\end{de}

\noindent
The Wronskian generalizes for suitably differentiable functions on $\Acal$ in the natural fashion.

\begin{de} {\bf Wronskian} of functions $y_1,y_2, \dots , y_m$ at least $(m-1)$ times differentiable at $z$ is given by:
\[ W(y_1,y_2,\dots , y_m; z) = det
\left[ 
	\begin{array}{cccc}
		y_1(z) & y_2(z) & \cdots & y_m(z)   \\
		y_1'(z) & y_2'(z) & \cdots & y_m'(z)   \\
\vdots &  \vdots & \cdots & \vdots  \\
		y_1^{(m-1)}(z) & y_2^{(m-1)}(z) & \cdots & y_m^{(m-1)}(z)   
	\end{array} 
\right].
\]  
\end{de}

\noindent
Notice that the Wronskian is formed by the determinant of a matrix of $\Acal$-elements for a given $z$. Fortunately, linear algebra over a commutative ring allows the usual theory of determinants. In particular, $\text{det}: \Acal^{ m \times m} \rightarrow \Acal$ and $M \in \Acal^{ m \times m}$ is invertible if and only if $\text{det}(M) \in \Acalx$. Furthermore, $Mx=0$ has nontrivial solutions if and only if $\text{det}(M) \in \Acalzd$. 

\begin{thm} 
Let $I \subseteq \Acal$ and suppose $y_1, \dots, y_m : I \rightarrow \Acal$ are at least $(m-1)$-times $\Acal$-differentiable. If $W(y_1,\dots , y_m; z) \in \Acalx$ for each $z \in I$ then $\{y_1, \dots , y_m\}$ is linearly independent on $I$. 
\end{thm}

\color{black}
\noindent
{\bf Proof:} Suppose for all $z \in I$ 
\begin{equation} \label{eqn:LIconditionI}
 c_1y_1(z)+c_2y_2(z)+ \cdots +c_my_m(z) = 0. 
 \end{equation}
Differentiate $(m-1)$ times to produce the following system of equations over $\Acal$:
\begin{equation} 
\underbrace{\left[ 
	\begin{array}{cccc}
		y_1(z) & y_2(z) & \cdots & y_m(z)   \\
		y_1'(z) & y_2'(z) & \cdots & y_m'(z)   \\
\vdots &  \vdots & \cdots & \vdots  \\
		y_1^{(m-1)}(z) & y_2^{(m-1)}(z) & \cdots & y_m^{(m-1)}(z)   
	\end{array} 
\right]}_{Y(z)} \left[ \begin{array}{l} c_1 \\ c_2 \\ \vdots \\ c_m \end{array} \right] = \left[ \begin{array}{l} 0 \\ 0 \\ \vdots \\ 0 \end{array} \right].
\end{equation}
Thus, Equation \ref{eqn:LIconditionI} has only the zero solution if and only if $\text{det}(Y(z)) \in \Acalx$ for each $z \in I$. But, this means $\{ y_1, \dots , y_m \}$ is linearly independent on $I$ if and only if $W(y_1, \dots, y_n ; z)$ is a unit for each $z \in I$. $\Box$ \\

\noindent
In practice, we are primarily interested in solution sets to linear $n$-th order $\Acal$-ODEs where the study of linear independence is greatly simplified by {\bf Abel's formula}. In particular, this formula forces the Wronskian of a full solution set to remain in either $\Acalx$ or $\Acalzd$ throughout the entirety of a connected subset. 

\begin{thm} \label{thm:abelsformula}
Suppose $a_o, a_1, \dots , a_n$ are continuous functions on the connected set $I \subseteq \Acal$ where $a_o(z) \in \Acalx$ for each $z \in I$.  If $y_1,y_2, \dots , y_n$ are solutions of $a_o y^{(n)}+ a_1 y^{(n-1)} + \cdots + a_{n-1} y' + a_ny=0$ then $ \displaystyle W(y_1, \dots, y_n;z) = C \, exp \left[ \int \frac{a_1}{a_o} \, dz \right] $ for each $z \in I$.
\end{thm}

\noindent
{\bf Proof:}  suppose $y_1,y_2, \dots , y_n$ are solutions on $I$ for $a_o y^{(n)}+ a_1 y^{(n-1)} + \cdots + a_{n-1} y' + a_ny=0$. Let $Y = [y_1,y_2,\dots , y_n]$ thus $Y' =  [y_1',y_2',\dots , y_n']$ and $Y^{(n-1)}= [y_1^{(n-1)},y_2^{(n-1)},\dots , y_n^{(n-1)}]$. The determinant which forms Wronskian is given by
\begin{equation}
 W = \sum_{i_1,i_2,\dots,i_n=1}^n \hspace{-0.1in} \eps_{i_1i_2...i_n} Y_{i_1}Y'_{i_2} \cdots Y^{(n-1)}_{i_n} 
\end{equation}
where $\eps_{i_1i_2...i_n}$ denoted the completely antisymmetric symbol where $\eps_{12\dots n}=1$. 
Apply the product rule for $n$-fold products on each summand in the above sum, 
\begin{equation}
 W' = \sum_{i_1,\dots,i_n=1}^n  \hspace{-0.1in} \eps_{i_1i_2...i_n}\bigl( Y_{i_1}'Y'_{i_2} \cdots Y^{(n-1)}_{i_n} 
+ Y_{i_1}Y''_{i_2}Y_{i_3}'' \cdots Y^{(n-1)}_{i_n} + \cdots +Y_{i_1}Y'_{i_2} \cdots Y^{(n-2)}_{i_{n-1}}Y^{(n)}_{i_n} \bigr).
\end{equation}
The term $Y_{i_1}'Y'_{i_2} \cdots Y^{(n-1)}_{i_n} = Y_{i_2}'Y'_{i_1} \cdots Y^{(n-1)}_{i_n}$ hence is  symmetric in the pair of indices $i_1,i_2$. Next, the term $Y_{i_1}Y''_{i_2}Y_{i_3}'' \cdots Y^{(n-1)}_{i_n}$ is symmetric in the pair of indices $i_2,i_3$. This pattern continues up to the term $Y_{i_1}Y'_{i_2} \cdots Y^{(n-1)}_{i_{n-2}}Y^{(n-2)}_{i_{n-1}}Y^{(n-1)}_{i_n}$ which is symmetric in the $i_{n-2},i_{n-1}$ indices. Thus all the terms vanish when contracted against the antisymmetric symbol. Only one term remains in calculation of $W'$:
\begin{equation} \label{eqn:wcancelling}
 W' = \sum_{i_1,i_2,\dots,i_n=1}^n \hspace{-0.1in} \eps_{i_1i_2...i_n} Y_{i_1}Y'_{i_2} \cdots Y^{(n-2)}_{i_{n-1}}Y^{(n)}_{i_n}  
 \end{equation}
Recall that $y_1,y_2, \dots , y_n$ are solutions of $a_o y^{(n)}+ a_1 y^{(n-1)} + \cdots + a_{n-1} y' + a_ny=0$ hence
\begin{equation}
 Y^{(n)} = -\frac{a_1}{a_o} Y^{(n-1)} - \cdots -\frac{a_{n-1}}{a_o}  Y' -\frac{a_n}{a_o} Y 
\end{equation}
Substitute this into Equation \ref{eqn:wcancelling},
\begin{align} 
W' &= \sum_{i_1,i_2,\dots,i_n=1}^n \hspace{-0.1in} \eps_{i_1i_2...i_n} Y_{i_1}Y'_{i_2} \cdots Y^{(n-2)}_{i_{n-1}}\biggl[ -\frac{a_1}{a_o} Y^{(n-1)} - \cdots -\frac{a_{n-1}}{a_o}  Y' -\frac{a_n}{a_o} Y \biggr]_{i_n}  \\ \notag 
&= \sum_{i_1,i_2,\dots,i_n=1}^n \hspace{-0.1in} \eps_{i_1i_2...i_n} \biggl( -\frac{a_1}{a_o} Y_{i_1}Y'_{i_2} \cdots Y^{(n-1)}_{i_n} - \cdots -\frac{a_{n-1}}{a_o}  Y_{i_1}Y'_{i_2} \cdots Y'_{i_n} -\frac{a_n}{a_o} Y_{i_1}Y'_{i_2} \cdots Y_{i_n} \biggr) \\ \notag 
&= -\frac{a_1}{a_o} \biggl( \, \sum_{i_1,i_2,\dots,i_n=1}^n  \hspace{-0.1in} \eps_{i_1i_2...i_n}   Y_{i_1}Y'_{i_2} \cdots Y^{(n-1)}_{i_n} \, \biggr) \qquad \star  \\ \notag 
&= -\frac{a_1}{a_o}W. 
\end{align}
The $\star$ step is based on the observation that the index pairs $i_1,i_n$ and $i_2, i_n$ etc... are symmetric in the line above it hence as they are summed against the completely antisymmetric symbol those terms vanish. Finally, we find $W' = -\frac{a_1}{a_o}W$ and conclude Abel's formula $\displaystyle W(y_1, \dots, y_n;z) = C \, exp \left[ -\int \frac{a_1}{a_o} \, dz \right]$ follows by integration since $I$ is connected\footnote{If $I$ was formed by several connected components then we could have different values for $C$ in different components.}.  $\Box$ \\

\noindent
In a field the only divisor of zero is zero itself hence we need only worry the Wronskian be zero in the ordinary theory. In $\Acal$-calculus we must also beware of nontrivial divisors of zero.

\begin{coro} \label{thm:abelsLI}
Suppose $a_o, a_1, \dots , a_n$ are continuous functions on the connected set $I \subseteq \Acal$ where $a_o(z) \in \Acalx$ for each $z \in I$.  Let $y_1,y_2, \dots , y_n$ be solutions of $a_o y^{(n)}+ a_1 y^{(n-1)} + \cdots + a_{n-1} y' + a_ny=0$. There exists $z_o \in I$ such that $W(y_1,y_2, \dots , y_n;z_o) \in \Acalx$ if and only if $\{y_1,y_2, \dots , y_n \}$ is linearly independent on $I$.  Likewise, there exists $z_o \in I$ such that $W(y_1,y_2, \dots , y_n;z_o) \in \Acalzd$ if and only if $\{y_1,y_2, \dots , y_n \}$ is linearly dependent on $I$.  
\end{coro}

\noindent
\color{black}
{\bf Proof:} Theorem \ref{thm:abelsformula} provides $ W(y_1, \dots, y_n;z) = C \, exp \left[ -\int \frac{a_1}{a_o} \, dz \right] $ for all $z \in I$. If there exists $z_o \in I$ such that $W(y_1,y_2, \dots , y_n;z_o) = Cexp \left[ -\int \frac{a_1}{a_o} \, dz \right]\bigg{|}_{z=z_o} \in \Acalzd$ then we find $C \in \Acalzd$ since the image of the exponential is in $\Acalx$. Likewise, if $W(y_1,y_2, \dots , y_n;z_o) \in \Acalx$ then $C \in \Acalx$. Thus the Wronskian of a solution set on a connected subset is either always a zero divisor or always a unit. $\Box$ 

\begin{de} 
Suppose $L[y]=f$ is an $n$-th order linear differential equation on connected $I \subseteq \Acal$. We say $S=\{ y_1,y_2, \dots, y_n \}$ is a {\bf fundamental solution set} of $L[y]=f$ if and only if $S$ is a linearly independent set of solutions to the homogeneous equation; $L[y_j]=0$ for $j=1,2,\dots n$. 
\end{de}

\noindent
Note the fundamental solution set of $L[y]=f \neq 0$  does not solve $L[y]=f$. We should mention the usual theory for nonhomogeneous differential equations is also naturally generalized to $\Acal$-calculus. We leave explicit discussion to a future work.

\begin{thm} \label{thm:fundsolutionsetexists}
If $L[y]=f$ is an $n$-th order linear differential equation with continuous coefficient functions on connected $I \subseteq \Acal$ then there exists a fundamental solution set $S=\{ y_1,y_2, \dots, y_n \}$ on $I$.  
\end{thm}

\noindent
\color{black}
{\bf Proof:} Apply Theorem \ref{thm:existenceIVP} $n$-times as to select $z_o \in I$ and unique solutions $y_1,\dots , y_n$ for which  $y_i^{(j)}(z_o) = \delta_{i,j-1}$ for $0 \leq i \leq n-1$ and $i=1,\dots , n$. Let the Wronskian at $z=z_o$ for the solution set $\{ y_1,y_2, \dots, y_n \}$ be $W(z)$ for the remainder of this proof:
\[ W(z_o) = det
\left[ 
	\begin{array}{cccc}
		y_1(z_o) & y_2(z_o) & \cdots & y_m(z_o)   \\
		y_1'(z_o) & y_2'(z_o) & \cdots & y_m'(z_o)   \\
\vdots &  \vdots & \cdots & \vdots  \\
		y_1^{(n-1)}(z_o) & y_2^{(n-1)}(z_o) & \cdots & y_n^{(n-1)}(z_o)   
	\end{array} 
\right]
= 
det
\left[ 
	\begin{array}{cccc}
		1 & 0 & \cdots & 0   \\
		0 & 1 & \cdots & 0   \\
\vdots &  \vdots & \cdots & \vdots  \\
		0 & 0 & \cdots & 1   
	\end{array} 
\right]=1.
\] 
Therefore, by Corollary \ref{thm:abelsLI} the solution set is linearly independent on $I$. $\Box$ 

\begin{thm}  \label{thm:homo}
If $L[y]=0$ is an $n$-th order linear differential equation with continuous coefficient functions and fundamental solution set $S=\{ y_1,y_2, \dots, y_n \}$ on connected $I \subseteq \Acal$. Then if $y$ solves $L[y]=0$ then there exist unique constants $c_1,c_2, \dots , c_n \in \Acal$ such that:
\[ y = c_1y_1+c_2y_2+ \cdots +c_ny_n. \] 
\end{thm}

\noindent
\color{black}
{\bf Proof:} Suppose $L[y]=0$. If $y = c_1y_1+c_2y_2+ \cdots c_ny_n$ and $z_o \in I$ then note  $y^{(j)}(z_o) = (c_1y_1+c_2y_2+ \cdots c_ny_n)^{(j)}(z_o)$ for $j=0,1,\dots , n-1$. That is, we must solve:
\begin{equation}
 \left[ \begin{array}{c} y(z_o) \\ y'(z_o) \\ \vdots \\ y^{(n-1)}(z_o) \end{array} \right] = 
\left[ 
	\begin{array}{cccc}
		y_1(z_o) & y_2(z_o) & \cdots & y_n(z_o)   \\
		y_1'(z_o) & y_2'(z_o) & \cdots & y_n'(z_o)   \\
\vdots &  \vdots & \cdots & \vdots  \\
		y_1^{(n-1)}(z_o) & y_2^{(n-1)}(z_o) & \cdots & y_n^{(n-1)}(z_o)   
	\end{array} 
\right]
 \left[ \begin{array}{l} c_1 \\ c_2 \\ \vdots \\ c_n \end{array} \right] 
\end{equation}
for $c_1,\dots, c_n \in \Acal$. Since $\{ y_1, \dots, y_n \}$ is a fundamental solution set we know the determinant of the coefficient matrix above is a unit (it is the Wronskian of $y_1,\dots , y_n$ at $z_o$) hence this system of equations has a unique solution.  $\Box$ \\

\section{first order differential equations over $\Acal$} \label{sec:firstorderAODEs}
The essential point which we illustrate here is that the usual methods given in the introductory course equally well apply to differential equations in $\Acal$. That said, we hope the examples illustrate the novel nature of such problems. 

\begin{exa}
Let $z,w$ denote hyperbolic variables. To solve $\frac{dw}{dz} = w^2 \sin(z)$ we note formally\footnote{in calculus over $\RN$, the chain rule allows us to separate variables, the same is true for calculus over $\Acal$.} $\frac{dw}{w^2} = \sin(z) dz$ hence integration yields $\frac{1}{w} = \cos(z)+c$ where $c \in \Hcal$. Thus,
\begin{equation}
w = \frac{1}{\cos (z)+c} .
\end{equation}
Let $z=x+jy$ and $c=a+bj$ where $x,y,a,b \in \RN$ and calculate:
\begin{equation}
w = \frac{a+\cos(x)\cos(y)+j(b-\sin(x)\sin(y))}{(a+\cos(x)\cos(y))^2-(b-\sin(x)\sin(y))^2}. 
\end{equation}
If we denote $w=u+jv$ where $u,v$ are real variables then the equation above reveals
\begin{align} 
u &= \frac{a+\cos(x)\cos(y)}{(a+\cos(x)\cos(y))^2-(b-\sin(x)\sin(y))^2}  \\ \notag
v &= \frac{b-\sin(x)\sin(y)}{(a+\cos(x)\cos(y))^2-(b-\sin(x)\sin(y))^2}. 
\end{align}
Our solution assumes $w=u+jv$ is an $\Hcal$-differentiable function of $zx+jy$ hence we implicitly impose the $\Hcal$-Cauchy Riemann equations $u_x=v_y$ and $u_y=v_x$. Moreover, since $\frac{dw}{dz} = \frac{\partial w}{\partial x}$ the initial $\Acal$-ODE $\frac{dw}{dz} = w^2 \sin(z)$ amounts to the real PDEs,
\begin{align} 
u_x &= (u^2+v^2)\sin x \cos y + 2uv \sin y \cos x \\ \notag
v_x &= 2uv\sin x \cos y + (u^2+v^2) \sin y \cos x.
\end{align}
\end{exa}

\begin{exa}
Consider the $3$-hyperbolic numbers $\mathcal{H}_3 = \RN \oplus j\RN\oplus j^2 \RN$ with variables $\zeta = x+jy+j^2z$ and $\eta = u+jv+j^2w$. The $\mathcal{H}_3$-ODE given by $\frac{d\eta}{d\zeta} = \eta$ has natural solution $\eta = ke^{\zeta}$ where $k = a+bj+cj^2$ for some $a,b,c \in \RN$. This solution implicitly solves:
\begin{equation} \label{eqn:systempdes}
\underbrace{u_x=u, \ v_x=v, \ w_x=w }_{ from \ \frac{d\eta}{d\zeta} = \eta} \ \ \& \ \ \underbrace{u_x=v_y=w_z, \ u_y=v_z=w_x, \ u_z=v_x=w_y}_{\mathcal{H}_3 \ CR-equations}. 
\end{equation}
The special functions of $\mathcal{H}_3$ are $\cosh_3, \sinh_{31}, \sinh_{32}$ where
\begin{equation}
e^{j \theta} = \cosh_3(\theta)+
j\sinh_{31}(\theta)+j^2\sinh_{32}(\theta). 
\end{equation}
Noting that $\cosh_3(j \theta) = \cosh_3(\theta)$ and $\sinh_{31}(j\theta) = j \sinh_{31}(\theta)$ and $\sinh_{32}(j \theta) = j^2 \sinh_{32}(\theta)$ as shown in \cite{bedell2} we derive the following component expansion for the exponential:
\begin{align} \label{eqn:nastyfla}
e^{\zeta} = e^{x+jy+j^2z} &= e^xe^{jy}e^{j^2z} \\ \notag
&= e^x\bigl[ \underbrace{\cosh_3(y)\cosh_3(z)+\sinh_{31}(y)\sinh_{31}(z) + +\sinh_{32}(y)\sinh_{32}(z)}_{g_1} \bigr]  \\ \notag
&\ \ \ \ + je^x\bigl[ \underbrace{\cosh_3(y)\sinh_{32}(z)+\sinh_{31}(y)\cosh_3(z) +\sinh_{32}(y)\sinh_{31}(z)}_{g_2} \bigr]  \\ \notag
&\ \ \ \ + j^2e^x\bigl[ \underbrace{\cosh_3(y)\sinh_{31}(z)+\sinh_{32}(y)\cosh_3(z) +\sinh_{31}(y)\sinh_{32}(z)}_{g_3} \bigr].  
\end{align}
Thus $e^{x+jy+j^2z} = e^x g_1+je^xg_2+j^2e^xg_3$. Hence,
\begin{align}
\eta = ke^{\zeta} &= e^x(a+bj+cj^2)( g_1+jg_2+j^2g_3) \\ \notag
&= e^x(ag_1+bg_3+cg_2)+j e^x(ag_2+bg_1+cg_3)+j^2e^x(ag_3+bg_2+cg_1).
\end{align}
Therefore, the system of PDEs given in Equation \ref{eqn:systempdes} finds solutions of the form:
\begin{equation}
u=e^x(ag_1+bg_3+cg_2), \ v= e^x(ag_2+bg_1+cg_3), \ w=e^x(ag_3+bg_2+cg_1).
\end{equation}
\end{exa}

\noindent
Even something as simple as an exponential solution in the algebra may have rather complicated real content. One ultimate goal of our work is to answer the {\it inverse problem}. In particular, we would like to learn when a given system of real PDEs can be recast as a problem of $\Acal$-ODEs for an appropriate choice of algebra. The present work does not seek to solve the inverse problem. Rather, we work towards gaining a deeper understanding of solution techniques for $\Acal$-ODEs. \\

\noindent
In the study of real first order ODEs one often begins with the study of separable, linear and exact differential equations. Since calculus over $\Acal$ has the same basic rules we naturally generalize the standard methods: 

\begin{thm}
Let $\zeta, \eta$ denote variables in $\Acal$ and suppose $f(\zeta),g(\eta)$ are continuous functions then $\frac{d\eta}{d\zeta} = f(\zeta)g(\eta)$ has solutions given implicitly by $\int \frac{1}{g(\eta)} d\eta =  \int f(\zeta) d\zeta$.
\end{thm}

\begin{thm}
Let $\zeta, \eta$ denote variables in $\Acal$ and suppose $P,Q$ are continuous functions then $\frac{d\eta}{d\zeta}+P\eta = Q$ has general solution $\eta = \frac{1}{I}\int IQ d\zeta$ where $I = \text{exp}( \int P d\zeta)$.
\end{thm}

\begin{thm}
Let $\zeta, \eta$ denote variables in $\Acal$ and suppose $M,N$ are $\Acal$-differentiable functions of $\zeta, \eta$. If there exists a function $F: \Acal \times \Acal \rightarrow \Acal$ for which $M = \frac{\partial F}{\partial \zeta}$ and $N = \frac{\partial F}{\partial \eta}$ then $Md\zeta+Nd\eta =0$ has solution $F(\zeta, \eta) = C$.
\end{thm}

\noindent
With the results above we can solve certain systems of real PDEs, even nonlinear, usually coupled, by elementary calculus in an appropriate algebra variable. 

\begin{exa}
The solution to $2\eta d\eta+2 \zeta d\zeta = 0$ in $\Hcal$ is implicitly given by the solution set of $\eta^2+\zeta^2=c$. In terms of real variables $x,y,u,v$ for which $\zeta = x+jy$, $\eta = u+jv$ and constants $a,b$ with $c=a+jb$ we find (as $j^2=1$)
\begin{equation}
x^2+y^2+u^2+v^2=a \ \ \& \ \ 2xy+2uv = b.
\end{equation}
The equations above solve $\frac{d\eta}{d\zeta} = \frac{-\zeta}{\eta}$ which amounts to  the following system of real PDEs,
\begin{equation}
\underbrace{u_x = \frac{yv-xu}{u^2-v^2}, \  v_x = \frac{xv-yu}{u^2-v^2} }_{from \ \frac{d\eta}{d\zeta} = \frac{-\zeta}{\eta}}, \ \ \& \ \ \underbrace{u_x=v_y, \ v_x=u_y}_{\Hcal \ CR-equations}.
\end{equation}
\end{exa}

\noindent
The choice of $\Acal = \Hcal$ is merely for illustration in the above example. We could replace $\Hcal$ with any algebra $\Acal$ and the underlying real PDEs could be much more complicated. That said, since we assume the existence of $\Acal$-derivatives there is a certain pairing of real coordinates which must be seen both in the solutions and the PDEs. Therefore, we expect that many systems of PDEs will not admit an interesting algebraic reformulation. 

\section{nondegenerate constant coefficient $\Acal$-ODEs} \label{sec:nondegencc}
A {\it constant coefficient} $\Acal$-ODE of $n$-th order has the form:
\begin{equation} \label{eqn:degeneratedefn}
a_n\frac{d^n \eta}{d\zeta^n}+ a_{n-1}\frac{d^{n-1} \eta}{d\zeta^{n-1}}
+ \cdots + 
a_1\frac{d\eta}{d\zeta}+a_0\eta = 0
\end{equation}
for {\it coefficients} $a_n,a_{n-1}, \dots , a_1, a_0 \in \Acal$ with $a_n \neq 0$. If $a_n \in \Acalx$ then we say the $\Acal$-ODE is {\bf nondegenerate} and if $a_n \in \Acalzd$ then the $\Acal$-ODE is {\bf degenerate}. We study the degenerate case in Section \ref{sec:degenerate}. We assume $a_n=1$ in most of what follows without loss of generality. \\

\noindent
Our discussion is divided into three main story arcs. First, we describe how to solve a constant coefficient $\Acal$-ODE via a given operator factorization in the ring $\Acal[ D]$ where $D = d/d\zeta$. Second, we study how isomorphism of algebras allows elegant solutions when the given algebra is known to be isomorphic to a direct product. Third, we introduce a novel generalization of the complexification technique. In particular, it is shown how the exponential on the natural extension algebra produces a fundamental solution set. We should emphasize from the outset, Theorems \ref{thm:fundsolutionsetexists} and \ref{thm:homo} show us that the general solution exists. We simply provide methods of calculation which reveal the explicit structure of that general solution.

\subsection{solution by operator technique}
In this section we assume $\Acal$ is a commutative algebra and denote $D = d/d\zeta$ for the operation of differentiation with respect to the algebra variable $\zeta$. We define $\Acal[D]$ to be the set of operators of the form
\begin{equation}
L = a_nD^n+a_{n-1}D^{n-1}+ \cdots + a_1D+a_0
\end{equation}
where $a_n,a_{n-1}, \dots , a_1,a_0 \in \Acal$. Observe $\Acal[D]$ forms an algebra with multiplication given by composition of operators. Commutativity of $\Acal$ and linearity of $D$ imply $\Acal[D]$ is also a commutative algebra.  Observe $\Acal[D]$ and $\Acal[z]$ are isomorphic via the map $\Psi(P(z))=P(D)$. In particular, this isomorphism allows us to pass factorizations in $\Acal[z]$ to corresponding factorizations in $\Acal[D]$:

\begin{thm} \label{thm:factorbysolution}
Suppose $D = d/d\zeta$ is the operator of differentiation with respect to the algebra variable $\zeta$. If $L \in \mathcal{A}[D]$ and there exists $\alpha \in \mathcal{A}$ such that $L[e^{\alpha \zeta}] = 0$ then there exists $L' \in \mathcal{A}[D]$ such that $L = (D-\alpha)L'$.
\end{thm}

\noindent
{\bf Proof:} Note $D^k[e^{\alpha \zeta}] = \alpha^k e^{\alpha \zeta}$ for each $k \in \NN$. Thus, if $L = P(D) = a_n D^n + \dots + a_1 D + a_0$, then $L[e^{\alpha \zeta}] = a_n \alpha^n e^{\alpha \zeta} + \dots + a_1 \alpha e^{\alpha \zeta} + a_0 e^{\alpha \zeta} = (a_n \alpha^n + \dots + a_1 \alpha + a_0)e^{\alpha \zeta}$.  Thus $L[e^{\alpha \zeta}]=0$ implies $a_n \alpha^n + \dots + a_1 \alpha + a_0 = 0$ thus $P(z) = a_n z^n + \dots + a_1 z + a_0 = (z-\alpha)g(z)$ for some $g(z) \in \Acal[z]$. Consequently,
\begin{equation}
\Psi (P(z)) = \psi ((z-\alpha)g(z)) = \psi(z-\alpha) \psi(g(z)) \ \ \Rightarrow \ \ P(D) = (D-\alpha)g(D)
\end{equation}
identify $L'=g(D) \in \Acal[D]$ hence $L =(D-\alpha)L'$ which completes the proof. $\Box$ \\

\noindent 
Thus, interestingly enough, the question of the existence of linear factors for a polynomial  in $\Acal[z]$ is tied to the existence of exponential solutions to the corresponding $\Acal$-ODE. Conversely, given a factorization of the operator defining an $\Acal$-ODE we obtain exponential solutions.

\begin{thm} 
Suppose $L \in \mathcal{A}[D]$ where $D = d/d\zeta$ defines an $\Acal$-ODE $L[\eta]=0$. If $L = (D-\alpha)L'$ where $\alpha \in \Acal$ then $\eta = e^{\alpha \zeta}$ is a solution of $L[\eta]=0$. \end{thm}

\noindent
{\bf Proof:} Observe $(D-\alpha)L' = L'(D-\alpha)$. Furthermore, $(D-\alpha)[e^{\alpha \zeta}] =  0$ hence $L[e^{\alpha \zeta}] = L'[(D-\alpha)[e^{\alpha \zeta}]]= L'[0]=0$. $\Box$ \\

\noindent
The result above naturally extends to multiple distinct factors. In particular, if $L =(D-\alpha_1)(D-\alpha_2)\cdots (D-\alpha_m)L'$ where $\alpha_1, \alpha_2, \dots , \alpha_m$ are distinct numbers in $\Acal$ then $\eta_1 = e^{\alpha_1 \zeta}, \eta_2 = e^{\alpha_2 \zeta}, \dots , \eta_m = e^{\alpha_m \zeta}$ solve $L[\eta]=0$. To treat repeated linear factors we note a standard lemma naturally transfers to $\Acal$-calculus:

\begin{lem}
Given a differential operator $p(D) \in \Acal[D]$, $\alpha \in \Acal$, and $f \in C(\Acal)$, then: $$ p(D)[e^{\alpha \zeta} f] = e^{\alpha \zeta} p(D+ \alpha)[f]. $$
\end{lem}

\noindent
{\bf Proof:} By the product rule, we have: $D[e^{\alpha \zeta} f] = \alpha e^{\alpha \zeta} f + e^{\alpha \zeta} D[f] = e^{\alpha \zeta}(D + \alpha)[f]$. It is then easy to see by induction that for all $k \geq 1$ we have $D^k[e^{\alpha \zeta} f] = e^{\alpha \zeta}(D+ \alpha)^k[f]$. Consider $p(D) = a_n D^n + \dots + a_1 D + a_0 \in \Acal[D]$ and calculate:
\begin{equation}
\begin{aligned}
 p(D)[e^{\alpha \zeta} f] &= (a_n D^n + \dots + a_1 D + a_0)[e^{\alpha \zeta} f] \\ &= a_n D^n[e^{\alpha \zeta} f] + \dots + a_1 D[e^{\alpha \zeta} f] + a_0(e^{\alpha \zeta} f) \\ &= e^{\alpha \zeta} a_n (D + \alpha)^n[f] + \dots + e^{\alpha \zeta} a_1 (D + \alpha)[f] + a_0 (e^{\alpha \zeta} f) \\ &= e^{\alpha \zeta} (a_n (D+\alpha)^n[f] + \dots + a_1 (D+\alpha)[f] + a_0 f) \\ &= e^{\alpha \zeta} p(D + \alpha)[f]
\end{aligned}
\end{equation}

\begin{thm} \label{thm:solnsetrepeatlinearfactor}
If $D= d/d\zeta$ and $L \in \Acal[D]$   factors as $L = (D - \alpha)^k L'$ for some differential operator $L' \in \Acal[D]$ and some $\alpha \in \Acal$, then $e^{\alpha \zeta}, \zeta e^{\alpha \zeta}, \dots, \zeta^{k-1} e^{\alpha \zeta} $ solve $L[ \eta ] = 0$.
\end{thm}

\noindent
{\bf Proof:} For all $i = 0, \dots, k-1$ we have: 
\begin{equation*}
\begin{aligned}
L[e^{\alpha \zeta} \zeta^i] = L'((D-\alpha)^k[e^{\alpha \zeta} \zeta^i]) = L'(e^{\alpha \zeta}D^k[\zeta^i]) = L'[0] = 0
\end{aligned}
\end{equation*}
Since $D^k[\zeta^i] = 0$ for all $i = 0, \dots, k-1$. Therefore, $e^{\alpha \zeta}, \zeta e^{\alpha \zeta}, \dots, \zeta^{k-1} e^{\alpha \zeta}$ are all solutions to the $\Acal$-ODE $L[\eta] = 0$. $\Box$ 

\begin{thm}
If $D= d/d\zeta$ and $L \in \Acal[D]$   factors as follows: 
$$ L = (D-\alpha_1)^{m_1}(D-\alpha_2)^{m_2}\cdots (D-\alpha_k)^{m_k} $$
for $\alpha_1, \alpha_2, \dots , \alpha_k \in \Acal$ such that $\alpha_i - \alpha_j \in \Acalx$ for all $i \neq j$ and $m_1,m_2, \dots , m_k \in \NN$ then the general solution to $L[\eta]=0$ is given by $ \displaystyle \eta = \sum_{i=1}^k \sum_{j=1}^{m_i}c_{ij}\zeta^{j-1}e^{\alpha_i \zeta}$ where $c_{ij} \in \Acal$.
\end{thm}

\noindent
The necessity of the condition $\alpha_i - \alpha_j \in \Acalx$ for distinct $i,j$ is illustrated by Example \ref{exa:nonuniquefactoring}.  \\

\noindent
{\bf Proof:} We apply Theorem \ref{thm:homo} to see it suffices to show $S=\{ \zeta^{j-1}e^{\alpha_i \zeta} \ | \ 1 \leq i \leq k, \ 1 \leq j \leq m_i \}$ is a fundamental solution set. It is already clear $S$ is a solution set by $k$-fold application of Theorem \ref{thm:solnsetrepeatlinearfactor}. Next, we show $S$ is linearly independent. Let $c_{ij} \in \Acal$ such that 
\begin{equation} \label{eqn:LIcondition}
\sum_{i=1}^k \sum_{j=1}^{m_i}c_{ij}\zeta^{j-1}e^{\alpha_i \zeta} = 0.
\end{equation} 
Construct $T_1 = L/(D-\alpha_1)$. Observe that $(D-\alpha_i)^me^{\alpha_j \zeta} = (\alpha_j - \alpha_i)^me^{\alpha_j \zeta}$ and $(D-\alpha_1)^{m_1-1}[\zeta^{m_1-1}e^{\alpha_1 \zeta}] = D^{m_1-1}[\zeta^{m_1-1}] = (m_1-1)!$. Consequently, if we operate on Equation \ref{eqn:LIcondition} by $T_1$ and evaluate at $\zeta=0$ then we obtain:
\begin{equation}
c_{1,m_1}(\alpha_2-\alpha_1)^{m_2}(\alpha_3-\alpha_1)^{m_3}\cdots (\alpha_k-\alpha_1)^{m_k}(m_1-1)! = 0.
\end{equation}
Thus $c_{1,m_1}=0$ as we assume, for distinct $i,j$, the difference between $\alpha_i-\alpha_j$ is a unit. Next, repeat the calculation for $L/(D-\alpha_1)^2$ to obtain $c_{1,m_1-1}=0$. Similar calculations will show $c_{1,j}=0$ for $j=1,\dots , m_1$. Continuing in this fashion for $\alpha_2, \dots , \alpha_k$ we find $c_{ij}=0$ for $1 \leq j \leq m_i$ where $1 \leq i \leq k$. Therefore, $S$ is linearly independent. $\Box$ \\

\noindent
In summary, if $P(D) \in \Acal[D]$ is the product of linear factors in $\Acal[D]$ whose corresponding zeros differences are units then we have a clear path to construct a solution of $P(D)[\eta]=0$. Of course, in $\Acal=\RN$ or $\CN$ every nonzero number is a unit hence the condition on $\alpha_1, \dots , \alpha_k$ simply reduces to $\alpha_i \neq \alpha_j$ for distinct $i,j$. The subtlety here is tied to the fact that $\Acal[D]$ is not generally a unique factorization domain when $\Acal$ is not a field. There exist distinct $\alpha_1,\dots , \alpha_k \in \Acal$ whose exponentials $e^{\alpha_1 \zeta}, \dots , e^{\alpha_k \zeta}$ have linear dependence:

\begin{exa} \label{exa:nonuniquefactoring}
Consider the hyperbolic numbers $\Hcal$ with variable $z=x+jy$. Observe
$D^2-1$ may be factored as either $(D-1)(D+1)$, or $(D-j)(D+j)$ over $\mathcal{H}$, and hence we find $e^z, e^{-z}$ as well as $e^{j z}, e^{-j z}$ as seemingly distinct exponential solutions to $w''-w=0$. However,  the independence of these solutions is  illusory. A short calculation reveals:
\begin{equation}
e^{j z} = \frac{1}{2}(1+j)e^z + \frac{1}{2}(1-j)e^{-z},  \ \ \& \ \ 
e^{-j z} = \frac{1}{2}(1+j)e^{-z} + \frac{1}{2}(1-j)e^{z}.
\end{equation}
Algebraically, these equations are very interesting, they express a linear dependence amongst zero divisors $1 \pm j$ and units. Of course, we should expect this result. Recall,
Theorem \ref{thm:homo} indicates that an $n$-th order $\Acal$-ODE can have at most $n$ linearly independent solutions. 
\end{exa}

\noindent
We must take care to find hidden dependence between seemingly distinct factors. If $(D-\alpha)$ and $(D-\beta)$ have $\alpha-\beta$ is a zero-divisor then these factors are not independent. 

\begin{exa}
Consider the differential equation $(D-1)(D+j)(D-j)[w]=0$ in $\Hcal$ we have $1+j$ and $1-j$ are zero divisors. Note:
\begin{equation}
(D-1)(D+j)(D-j)[w]=0 \ \ \Rightarrow \ \ (D-1)(D^2-1)[w] = 0 \ \ \Rightarrow \ \ (D-1)^2(D+1)[w] =0.
\end{equation}
Thus $w = c_1e^z+c_2ze^z+c_3e^{-z}$ is the general solution. 
\end{exa}

\noindent
The calculation in the example above is fortunate. In contrast, zero divisors cannot be hidden completely in the calculation below:

\begin{exa} \label{exa:hyperbolicwithzddiffroot}
Consider the differential equation $(D-1)(D-j)[w]=0$ in $\Hcal$ we have $1+j$ is a zero divisor. 
It is simple to check that $e^z$ and $e^{jz}$ are solutions. However,
\begin{equation}
\text{det} \left[ \begin{array}{cc}
e^z & e^{jz} \\
e^z & je^{jz}
\end{array} \right] = (j-1)e^{(1+j)z} = (j-1)\left(1+(1+j)z+ \frac{1}{2}(1+j)^2z^2+ \cdots \right) = j-1
\end{equation}
indicates linear dependence. In fact, $(1-j)e^z+(1+j)e^{jz}=0$ provides the explicit dependence. We rely on substitution to derive the solution. Let $\eta = (D-j)[w]$ thus solve $(D-1)[\eta]=0$ to obtain $\eta = c_1e^z$. Hence,
\begin{equation} \label{eqn:intfacttechnique}
\frac{dw}{dz}-jw = c_1e^z \ \ \Rightarrow \ \ e^{-jz}\frac{dw}{dz}-je^{-jz}w = c_1e^ze^{-jz} \ \ \Rightarrow \ \  \frac{d}{dz} \left( e^{-jz}w \right) =  c_1e^{(1-j)z}.
\end{equation}
Since $1-j \in \Acalzd$ we use series techniques to derive the explicit solution,
 \begin{equation}
\int e^{(1-j)z} dz = \int \left[1+\sum_{n=1}^{\infty} \frac{1}{n!} (1-j)^nz^n \right] dz = c + z+
\sum_{n=1}^{\infty} \frac{1-j}{(n-1)!(n+1)} z^{n+1}
\end{equation}
where we used the identity $(1-j)^n = n(1-j)$ for $n \in \NN$. Integrating Equation \ref{eqn:intfacttechnique} yields:
 \begin{equation} \label{eqn:generalsol}
w = c_2e^{jz}+c_1e^{jz}\left( z+
\sum_{n=1}^{\infty} \frac{1-j}{(n-1)!(n+1)} z^{n+1}\right).
\end{equation}
Setting $w_1(z) = e^{jz}$ and $w_2(z) = e^{jz} \left(z +
\sum_{n=1}^{\infty} \frac{1-j}{(n-1)!(n+1)} z^{n+1}\right)$ have $w_1(0)=0$ and $w_1(0)=j$ whereas $w_2(0)=0$ and $w_2'(0)=j$. It follows the Wronskian of $w_1,w_2$ evaluated at $z=0$ yields $j \in \Hcal^{\times}$ hence $\{ w_1(z), w_2(z) \}$ form a fundamental solution set and Equation \ref{eqn:generalsol} expresses the general solution to $(D-1)(D-j)[w]=0$. 
\end{exa}

\noindent
Equation \ref{eqn:generalsol} shows that roots which differ by a zero divisor behave somewhat like repeated roots. It is interesting that the term which is truly unfamilar to the student of standard ODEs has a manifest zero divisor. Novel terms are often attached to zero divisors. \\

\noindent
Finally, a word of caution, there are many linear differential operators which do not factor over a given algebra. For example, $D^2+1$ over $\RN$. Or, some may only partially split like $D^3-1=(D-j)(D^2 + jD + j^2)$ over $\RN \oplus j \RN \oplus j^2 \RN$ where $j^3=1$. To deal with irreducible operators of degree two or higher we use the {\it extension technique} given in Section \ref{sec:stupidsolutions}. 

\subsection{solution via isomorphism}
Suppose $\Acal$ is formed by the direct product of algebras $\Acal_1, \Acal_2, \dots , \Acal_k$ whose multiplications are denoted by juxtaposition. In particular, $\Acal = \Acal_1 \times \Acal_2 \times \cdots \times \Acal_k$ and for $x,y \in \Acal$:
\begin{equation}
xy = (x_1, x_2, \dots , x_k)(y_1,y_2, \dots , y_k) = (x_1y_1, x_2y_2, \dots , x_ky_k).
\end{equation}
If we denote $\mathds{1}_j$ for the unity in $\Acal_j$ for $j=1,2,\dots , k$ then the unity $\mathds{1}$ of $\Acal$ is precisely $\mathds{1} = ( \mathds{1}_1, \mathds{1}_2, \dots , \mathds{1}_k )$. Suppose bases $\beta_1, \beta_2, \dots , \beta_k$ are given such that 
\begin{equation}
\beta_{j} = \{ v_{1j}, v_{2j}, \dots , v_{n_j,j} \}
\end{equation}
with $v_{1j} = \mathds{1}_{j}$  and $n_j = \text{dim}(\Acal_j)$ for $j=1,2, \dots, k$. Introduce notation
\begin{equation}
e_1 = (\mathds{1}_1, 0 , \dots , 0), \ 
e_2 = (0, \mathds{1}_2 , \dots , 0), \ \dots, \ 
e_k = (0, 0 , \dots , \mathds{1}_k)
\end{equation}
so that $\beta = e_1\beta_1 \cup \cdots \cup e_k\beta_k$ given the natural order provides a basis for $\Acal$ whose dimension $n = n_1+n_2+ \cdots + n_k$. Furthermore, the regular representations of $\Acal$ naturally connect to regular representations of the component algebras:
\begin{prop} \label{prop:blockdiagrep}
If $\Acal_j$ has basis $\beta_j$ for $j=1,2, \dots , k$ and $\Acal = \Acal_1 \times \Acal_2 \times \cdots \Acal_k$ is given basis $\beta = e_1\beta_1 \cup e_2\beta_2 \cup \cdots \cup e_k\beta_k$  then for each $z= (z_1,z_2 \dots , z_k) \in \Acal$ we find
$$ \MM_{\beta}(z) = \MM_{\beta_1}(z_1) \oplus \MM_{\beta_2}(z_2) \oplus \cdots \oplus \MM_{\beta_1}(z_k). $$ 
\end{prop}

\noindent
Let $\zeta_1, \zeta_2, \dots , \zeta_k$ denote variables in $\Acal_1, \Acal_2, \dots , \Acal_k$ respective. If $\zeta \in \Acal$ then we may write:
\begin{equation}
\zeta = e_1 \zeta_1+ e_2\zeta_2+ \cdots + e_k \zeta_k.
\end{equation}
Each algebra variable $\zeta_j$ can be expanded in terms of its real substructure with respect to basis $\beta_j  = \{ v_{1j}, v_{2j}, \dots , v_{n_j,j} \}$ of $\Acal_j$
\begin{equation}
\zeta_j = x_{1j}v_{1j}+x_{2j}v_{2j}+ \cdots + x_{n_j j}v_{n_j j}.
\end{equation}
It follows that $f: \Acal \rightarrow \Acal$ can be understood as a function of the real variables $x_{ij}$ where $1 \leq i \leq n_j$ for $1 \leq j \leq k$. Let $f: \Acal \rightarrow \Acal$ have component functions $f_j: \Acal \rightarrow  \Acal_j$ where $f = f_1e_1+f_2e_2+ \cdots f_ke_k$. The Jacobian matrix of $f$ has the form:
\begin{equation}
J_f = \left[ 
 \left[\frac{ \partial f}{\partial x_{11}} \right]_{\beta} \bigg{|}  \cdots \bigg{|} \left[ \frac{ \partial f}{\partial x_{n_1,1}} \right]_{\beta}\bigg{|} \cdots \bigg{|} \left[\frac{ \partial f}{\partial x_{1k}} \right]_{\beta} \bigg{|}  \cdots \bigg{|} \left[ \frac{ \partial f}{\partial x_{n_k,k}} \right]_{\beta} \right].
\end{equation}
If we suppose $f: \Acal \rightarrow \Acal$ is $\Acal$-differentiable at $p$ then we have $J_f(p)$ is in the regular representation of $\Acal$. Proposition \ref{prop:blockdiagrep} indicates that $J_f(p)$ is a block-diagonal matrix. Omitting the coordinate maps and use block notation we have:
\begin{equation}
J_f = \left[ \begin{array}{c|c|c|c}
\partial_{i_1,1} f_1 & \partial_{i_2,2} f_1 & \cdots & \partial_{i_k,n_k} f_1 \\ \hline
\partial_{i_1,1} f_2 & \partial_{i_2,2} f_2 & \cdots & \partial_{i_k,n_k} f_2 \\ \hline
\vdots & \vdots & \cdots & \vdots \\ \hline
\partial_{i_1,1} f_k & \partial_{i_2,2} f_k & \cdots & \partial_{i_k,n_k} f_k
\end{array}\right] =
\left[ \begin{array}{c|c|c|c}
\partial_{i_1,1} f_1 & 0 & \cdots & 0 \\ \hline
0 & \partial_{i_2,2} f_2 & \cdots & 0 \\ \hline
\vdots & \vdots & \cdots & \vdots \\ \hline
0 & 0 & \cdots & \partial_{i_k,n_k} f_k
\end{array}\right]
\end{equation}
where $1 \leq i_j \leq n_j$ for $j=1,2,\dots , k$. Observe, for $j=1,2,\dots , k$, $f_j$ is only a function of $x_{i_jj}$ for $i_j=1,2, \dots , n_j$. Proposition \ref{prop:blockdiagrep} indicates $[\partial_{i_j,j} f_j]$ is in the regular representation of $\Acal_j$ with respect to $\beta_j$. Thus, abusing notation slightly as to view $f_j: \Acal_j \rightarrow \Acal_j$, we find $d_{p_j} f_j  \in \mathcal{R}_{\Acal_j}$ for $j=1,2, \dots , k$. Hence, an $\Acal$-differentiable function $f: \Acal \rightarrow \Acal$ has the form:
\begin{equation}
f( \zeta_1, \zeta_2, \dots , \zeta_k) = e_1f_1  (\zeta_1)+e_2f_2(\zeta_2)+ \cdots +e_k f_k(\zeta_k).
\end{equation}
Since $\mathds{1} = e_1+e_2+\cdots + e_k$ and $d_pf_j(e_i) = 0$ for $i \neq j$ we find 
\begin{align}
 d_pf( \mathds{1} ) &=
e_1d_pf_1( \mathds{1} ) +
e_2d_pf_2( \mathds{1} ) + \cdots +
e_kd_pf_k( \mathds{1} ) \\ \notag
&= 
e_1d_pf_1( e_1 ) +
e_2d_pf_2(e_2 ) + \cdots +
e_kd_pf_k( e_k ) \\ \notag
&= 
e_1d_{p_1}f_1( \mathds{1}_1 ) +
e_2d_{p_2}f_2(\mathds{1}_2 ) + \cdots +
e_kd_{p_k}f_k( \mathds{1}_k ).
\end{align}
Therefore, as $d_pf( \mathds{1}) = \frac{df}{d\zeta}(p)$ and $d_{p_j}f_j( \mathds{1}_j) = \frac{df_j}{d\zeta_j}(p_j)$ we find 
\begin{equation}
\frac{df}{d\zeta}(p) = e_1\frac{df_1}{d\zeta_1}(p_1)+e_2\frac{df_2}{d\zeta_2}(p_2)+ \dots + e_k \frac{df_k}{d\zeta_k}(p_k). 
\end{equation}
Let us summarize our observations:
\begin{prop} \label{prop:operatorfactors}
If $f \in \FunA$ then we may view $f = (f_1, f_2, \dots , f_k)$ where $f_i \in \FunAj$. Furthermore, if we denote $D=\frac{d}{d\zeta}$ and $D_j = \frac{d}{d\zeta_j}$ then $D = (D_1,D_2, \dots , D_k)$ where 
$$ Df = (D_1,D_2, \dots , D_k)(f_1,\dots, f_k) = (D_1f_1,D_2f_2, \dots , D_kf_k). $$
\end{prop}

\noindent
Polynomials in the product algebra have a very simple structure. Observe
\begin{equation}
\zeta^n = (\zeta_1,\zeta_2, \cdots,\zeta_k)^n = (\zeta_1^n,\zeta_2^n, \cdots,\zeta_k^n).
\end{equation}
Therefore, given coefficients $a_i = (a_{i1},a_{i2}, \dots , a_{ik}) \in \Acal$ for $i=0,1,\dots , n$ the polynomial $p(\zeta) = a_n\zeta^n+ \cdots + a_1 \zeta + a_0$ has
\begin{equation} \label{eqn:polynomialofpolys}
p(\zeta) = \sum_{i=0}^n a_i \zeta^i =   \left(\sum_{i=0}^n a_{i1}\zeta_1^i, \ \sum_{i=0}^n a_{i2}\zeta_2^i, \ \dots , \ \sum_{i=0}^n a_{ik}\zeta_k^i \right)
\end{equation}
Consequently, applying Proposition \ref{prop:operatorfactors} to the polynomial above:
\begin{equation} \label{eqn:polynomialOps}
p(D) = \sum_{i=0}^n a_i D^i =   \left(\sum_{i=0}^n a_{i1}D_1^i, \ \sum_{i=0}^n a_{i2}D_2^i, \ \dots , \ \sum_{i=0}^n a_{ik}D_k^i \right).
\end{equation}

\noindent
This algebra allows us to solve an $n$-th order constant coefficient ODE on $\Acal= \Acal_1 \times \cdots \times \Acal_k$ by solving $k$-related constant coefficient $n$-th order  ODEs for $\Acal_1, \Acal_2, \dots , \Acal_k$.

\begin{thm} \label{thm:productODEs}
Suppose $p(D) = (p_1(D_1),p_2(D_2), \dots ,p_k(D_k))$ where $p(D) \in \Acal[D]$ and $p_j(D_j) \in \Acal_j[D_j]$ for $j=1,2,\dots , k$. Then $w=(w_1,w_2, \dots , w_k)$ is a solution to $p(D)[w]=0$ if and only if $p_j(D_j)[w_j]=0$ for $j=1,2,\dots , k$.
\end{thm}

\noindent
{\bf Proof:} by Proposition \ref{prop:operatorfactors} and $w=(w_1,w_2, \dots , w_k)$ then:
\begin{equation}
p(D)[w] = (p_1(D_1)[w_1],p_2(D_2)[w_2], \dots ,p_k(D_k)[w_k]).
\end{equation}
Thus $p(D)[w]=0$ if and only if $p_j(D_j)[w_j]=0$ for $j=1,2,\dots , k$. $\Box$ \\

\noindent
Usually the given algebra is not manifestly a direct product hence we have to filter the calculus through connecting isomorphisms. Let us describe how this works. The main tool is the chain rule which links derivatives in isomorphic algebras. If $\eta = \Psi^{-1} \comp w \comp \Psi$ where $\Psi: \Acal \rightarrow \Bcal$ is an isomorphism and $w \in \FunB$ then $\eta \in \FunA$ and $\eta^{(j)} =  \Psi^{-1} \comp w^{(j)} \comp \Psi$ where $\eta^{(j)}$ denotes the $j$-fold $\Acal$-derivative of $\eta$ and $w^{(j)}$ denotes the $j$-fold $\Bcal$-derivative of $w$ for $j=1,2,\dots$. Moreover, $\Psi \comp\eta^{(j)} =  w^{(j)} \comp \Psi$. Suppose there exist $a_0,a_1, \dots , a_n \in \Acal$ and $\eta \in \FunA$ such that
\begin{equation} \label{eqn:tokenAODE}
a_n{\eta}^{(n)}+\cdots +a_1{\eta}'+a_0\eta =0.
\end{equation}
Since $\eta^{(j)} =  \Psi^{-1} \comp w^{(j)} \comp \Psi$ we find
\begin{equation} \label{eqn:tokenAODEpartII}
a_n(\Psi^{-1} \comp w^{(n)} \comp \Psi)+\cdots +a_1(\Psi^{-1} \comp w' \comp \Psi)+a_0 (\Psi^{-1} \comp w \comp \Psi) =0.
\end{equation}
The result below follows naturally:

\begin{thm} \label{thm:mapofDEs}
We find $\eta \in \FunA$ solves $a_n{\eta}^{(n)}+\cdots +a_1{\eta}'+a_0\eta=0$ if and only if $w = \Psi \comp \eta \comp \Psi^{-1} \in \FunB$ solves $\Psi(a_n)w^{(j)} +\cdots  +\Psi(a_1)w' +\Psi(a_0)w =0
$.
\end{thm}

\noindent
One consequence of the Wedderburn-Artin Theorem is that any commutative semisimple associative algebra $\Acal \approxeq \RN^m \times \CN^k$. In particular, in \cite{bedell1} arguments are given to explain the following isomorphisms\footnote{$\Hcal_n = \RN \oplus j \RN \oplus \cdots \oplus j^{n-1} \RN$ where $j^n=1$ and $\mathcal{C}_n = \RN \oplus i \RN \oplus \cdots \oplus i^{n-1} \RN$ where $i^n=-1$.} 
\begin{equation}
(1.) \ \mathcal{H}_{2k} \cong \mathbb{R}^2 \times \mathbb{C}^{k-1} \qquad
(2.) \ \mathcal{C}_{2k} \cong \mathbb{C}^k \  \qquad
(3.) \ \mathcal{H}_{2k + 1} \cong \mathcal{C}_{2k + 1} \cong \mathbb{R} \times \mathbb{C}^k.
\end{equation}
Consequently, combining Theorems \ref{thm:productODEs} and \ref{thm:mapofDEs} we may solve a constant-coefficient $\Acal$-ODE for a commutative semisimple algebra by solving $m$-real and $k$-complex ODEs. 

\begin{exa} \label{exa:firstorderzerodivODE}
Consider $\Hcal$ and $\RN \times \RN$.  The isomorphism $\Psi( x+jy) = (x+y,x-y)$ maps $\Hcal$ to $\RN \times \RN$ has inverse $\Psi^{-1}(a,b) = \frac{1}{2}[(a+b)+j(a-b) ]$. Consider $\eta'-(1+j)\eta =0$ this maps to $w'-(2,0)w=0$ in $\RN \times \RN$. Denoting $s,t$ as the coordinates in $\RN \times \RN$ we find solution $w(s,t) = (c_1e^{2s}, c_2)$. Hence,
\begin{align}
\eta (x+jy) &= \Psi^{-1}(w(x+y,x-y)) \\ \notag
&= \Psi^{-1}(c_1e^{2(x+y)},c_2) \\ \notag
&= \frac{1}{2}[(c_1e^{2(x+y)}+c_2)+j(c_1e^{2(x+y)}-c_2) ] \\ \notag
&= \frac{c_1}{2}(1+j)e^{2(x+y)} + \frac{c_2}{2}(1-j).
\end{align} 
solves $\eta'-(1+j)\eta= 0$. Of course, we can also solve directly in $\Hcal$
\begin{equation}
\eta(x+jy) = Be^{(1+j)(x+jy)} 
\end{equation}
and a short calculation shows these solutions are equivalent.
\end{exa}

\begin{exa}
If $\zeta \in \mathcal{H}_4$ then $\zeta = t+jx+j^2y+j^3z$ for $t,x,y,z \in \RN$. Consider the $\mathcal{H}_4$ ODE given by $\eta''+j^2\eta =0$. Under the isomorphism $\Psi: \mathcal{H}_4 \rightarrow \RN^2 \times \CN $ given by $\Psi(j) = (1,-1,i)$ which implies for $t+xj+yj^2+zj^3 \in \mathcal{H}_4$:
\begin{equation}
\Psi( t+xj+yj^2+zj^3) = (t+x+y+z, \ t-x+y-z, \ t+ix-y-iz).
\end{equation}
Note $\Psi(j^2) = (\Psi(j))^2 = (1,-1,i)^2 = (1,1,-1)$ thus study the related $\RN^2 \times \CN$ ODE
\begin{equation}
w''+(1,1,-1)w = 0 \ \ \Rightarrow \ \ w_1''+w_1=0, \ w_2''+w_2=0, \ w_3''-w_3=0.
\end{equation}
Here $w_1' = \frac{dw_1}{dx_1}$ and $w_2' = \frac{dw_2}{dx_2}$ denote real derivatives and $w_3' = \frac{dw_3}{dx_3}$ is a complex derivative. Hence, by the usual arguments, we find general solutions for $w_1,w_2,w_3$ and hence:
\begin{equation}
w(x_1,x_2,x_3) = \left( c_1\cos( x_1)+ c_2\sin(x_1) , \  c_3\cos( x_2)+ c_4\sin(x_2) ,  \ c_5e^{x_3}+c_6e^{-x_3} \right)
\end{equation}
where $c_1,c_2,c_3,c_4 \in \RN$ and $c_5,c_6 \in \CN$. It follows we have solution $\eta = \Psi^{-1} \comp w \comp \Psi$ to solve $\eta'' + j^2 \eta = 0$. The formula for $\eta$ as a function of $t,x,y,z$ is lengthy so we omit it. That said, it may be interesting to relate the solution obtained via the isomorphism to the natural $\mathcal{H}_4$ solution $\eta = B_1\cos(j\zeta)+B_2\sin(j \zeta)$ where $B_1,B_2 \in \mathcal{H}_4$.
\end{exa}

\subsection{solution by extension} \label{sec:stupidsolutions}

\noindent
In ordinary ODEs (i.e. the usual theory of ODEs in $\mathbb{R}$), a standard technique is to factor the differential operator not over $\mathbb{R}[D]$, but $\mathbb{C}[D]$. For example, $D^2+1$ factors as $(D-i)(D+i)$ over $\mathbb{C}[D]$. Then, since $\mathbb{C}$ is a field then, the factored $(D-i)(D+i)[w] = 0$ provides $e^{iz}$ and $e^{-iz}$ as solutions to our original real ODE reinterpreted as a complex ODE. However, the more interesting aspect of this is that the real component functions of $e^{ix}$, namely, $\cos(x)$ and $\sin(x)$ in fact solve the original real ODE ! \\

\noindent
We wish to generalize this technique to arbitrary associative algebras, proving that the component functions of our generalization of the complex exponential give $n$ linearly independent solutions to the given differential equation, but we first must discuss some preliminaries.

\begin{de}
An extension algebra $\mathcal{A}'$ of $\mathcal{A}$ is an algebra which contains $\mathcal{A}$ as an isomorphic subalgebra. Furthermore, if $\mathcal{A}' = e_1\Acal \oplus e_2\Acal\oplus \cdots \oplus e_n \Acal$ where $\{ e_1, \dots , e_n \}$ is an $\Acal$-module basis for $\Acal'$. Suppose $U \subseteq \Acal$ is an open set. If $f: U \rightarrow \Acal'$ and $f = e_1f_1+e_2f_2+ \cdots +e_nf_n$ and we say $f_i: U \rightarrow \Acal$ for $i=1,2,\dots , n$ are the {\bf component functions}. We say $f: U \rightarrow \Acal'$ is {\bf $\Acal$-differentiable at $p \in U$} if each component function of $f$ is $\Acal$-differentiable at $p$. Furthermore, if $f$ is $\Acal$-differentiable for each $p \in U$ then $f$ is {\bf $\Acal$-differentiable on $U$}. Finally, if $f$ is differentiable on $U$ then we define the derivative function $\frac{df}{d\zeta}: U \rightarrow \Acal$ by:
$$ \frac{df}{d\zeta} = e_1\frac{df_1}{d\zeta}+e_2\frac{df_2}{d\zeta}+ \cdots +e_n\frac{df_n}{d\zeta}.$$
\end{de}

\noindent
It is helpful to appreciate how $\Acal$-differentiability relates to right-$\Acal$-linearity of the differential. Notice, if $f_i: \Acal \rightarrow \Acal$ is $\Acal$-differentiable at $p$ then $d_pf_i \in \EndA$ which means $d_pf_i(xy) = d_pf_i(x)y$ for all $x,y \in \Acal$. Thus, if $f = e_1f_1 + \cdots +e_nf_n: \Acal \rightarrow \Acal'$ is $\Acal$-differentiable at $p$ then for $x,y \in \Acal$,
\begin{align} \label{eqn:rightAofAprimemap}
(d_pf)(xy) &= e_1d_pf_1(xy) + \cdots + e_nd_pf_n(xy) \\ \notag
&= e_1d_pf_1(x)y + \cdots + e_nd_pf_n(x)y \\ \notag
&= \left(e_1d_pf_1(x) + \cdots + e_nd_pf_n(x) \right)y \\ \notag
&= d_pf(x)y.
\end{align}
Recall $\Acal$-differentiability of $f_i$ at $p$ indicates $f_i'(p) = d_pf_i(1)$ for $i=1,\dots , n$. Thus:
\begin{equation} \label{eqn:diffioneforprimemap}
f'(p) = e_1f_1'(p)+ \cdots + e_nf_n'(p) = 
e_1d_pf_1(1)+ \cdots + e_nd_pf_n(1) = d_pf(1).
\end{equation}

\noindent
The identities above in conjunction with the usual chain rule of multivariate real analysis yield the following mixed chain rule:\footnote{this is a generalization of the mixed chain rule from complex analysis: if $z \mapsto g(z)$ is complex differentiable and $t \mapsto f(t)$ is real differentiable then $g \comp f$ is real differentiable and $\frac{d}{dt}g(f(t)) = \frac{dg}{dz}(f(t))\frac{df}{dt}$}

\begin{thm} \label{thm:mixedchainrule}
Let $\Acal$ be an $m$-dimensional commutative, associative, unital real algebra. Let $\Acal'$ be a rank $n$ extension algebra of $\Acal$. Suppose $f: \Acal \rightarrow \Acal'$ is $\Acal$-differentiable at $p$ and $g: \Acal' \rightarrow \Acal'$ is $\Acal'$-differentiable at $f(p)$. Then $g \comp f$ is $\Acal$-differentiable at $p$ and
$$ (g \comp f)'(p) = g'(f(p))f'(p). $$
\end{thm}

\noindent
If we use $\zeta$ to denote $\Acal$-derivatives and $\zeta'$ to denote $\Acal'$-derivatives then the chain rule given above is written as:
\begin{equation}
\frac{d}{d\zeta}\bigg{|}_{\zeta=p}\left[g(f(\zeta)) \right] =
\frac{dg}{d\zeta'}(f(p)) \frac{df}{d\zeta}(p)
\end{equation}
where we should appreciate the multiplication of $\frac{df}{d\zeta}(p)$ and $\frac{dg}{d\zeta'}(f(p))$ is accomplished in $\Acal'$. \\

\noindent
{\bf Proof:} Let $\Acal'$ have $\Acal$-module basis $\{ e_1, \dots , e_n \}$. If $g$ has component functions $g_1,\dots , g_n$ then the $i$-th component function of the composite $g \comp f$ is simply $g_i \comp f: \Acal \rightarrow \Acal$. Note $g \comp f: \Acal \rightarrow \Acal'$ is real differentiable at $p$ as it is the composite of real differentiable maps. Therefore, by the usual chain rule of advanced calculus, we have $d_p(g \comp f) = d_{f(p)}g \comp d_pf$. Thus, for $x,y \in \Acal$, combining the chain rule and Equation \ref{eqn:rightAofAprimemap}
\begin{equation}
d_p(g_i \comp f)(xy) = d_{f(p)}g_i (d_pf(xy)) = d_{f(p)}g_i (d_pf(x)y).
\end{equation}
But, $y \in \Acal \subset \Acal'$ hence by $\Acal'$-differentiability of $g$ at $f(p)$ we find $d_{f(p)}g_i (d_pf(x)y) = d_{f(p)}g_i (d_pf(x))y$ which shows $d_p(g_i \comp f) \in \EndA$. Thus $g_i \comp f$ is $\Acal$-differentiable at $p$ for $i=1,\dots , n$ and this establishes $g \comp f: \Acal \rightarrow \Acal'$ is $\Acal$-differentiable. Hence,
\begin{equation}
 d_p(g \comp f)(1) = d_{f(p)}g(d_pf(1))  = d_{f(p)}g(1)d_pf(1) = g'(f(p))f'(p). 
\end{equation}
Therefore, $(g \comp f)'(p) = g'(f(p))f'(p)$. $\Box$ \\

\noindent
The exponential on $\Acal'$ is $\Acal'$-differentiable and, given a constant $\alpha \in \Acal'$, the map $\zeta \mapsto \alpha \zeta$ is $\Acal$-differentiable map from $\Acal$ to $\Acal'$. Applying Theorem \ref{thm:mixedchainrule} we find:

\begin{coro} \label{coro:exponAprime}
Suppose $\Acal'$ is a rank $n$ extension algebra of $\Acal$ then for $\alpha \in \Acal'$ 
$$ \frac{d}{d\zeta} \text{exp}( \alpha \zeta) = \alpha \text{exp}( \alpha \zeta).$$
\end{coro}

\noindent
The formula above is very interesting when we apply it to the extension algebra which is characteristic to a given constant coefficient $\Acal$-ODE. 

\begin{de}
Let $p(x) \in \Acal[x]$ be a monic $n^{th}$-order polynomial. The constant coefficient $n$-th order $\Acal$-ODE $p(D)[\eta]=0$ has  {\bf characteristic extension algebra} defined as follows:
\begin{equation}
\mathcal{A}' = \frac{\Acal[x]}{\langle p(x) \rangle } = \{ x_0+x_1k+\cdots +x_nk^{n-1} \ | \ x_0,x_1,\dots, x_n \in \Acal, \ k = x+\langle p(x) \rangle \}
\end{equation}
Furthermore,  $f: \Acal \rightarrow \Acal'$ defined by $f(\zeta) = \text{exp}( k \zeta)$ has {\it component functions} $f_1,\dots , f_n: \Acal \rightarrow \Acal$ defined implicitly by:
\begin{equation}
\text{exp}( k \zeta) = f_1( \zeta)+f_2(\zeta)k+ \cdots + k^{n-1}f_n(\zeta).
\end{equation}
We say $f_1,f_2, \dots , f_n$ as above are the {\bf special functions} of exponential on $\Acal'$.
\end{de}

\noindent
We say $\{ 1,k, \dots , k^{n-1} \}$ forms the {\bf standard $\Acal$-module basis} for $\Acal'$. Linear independence of $\{ 1,k, \dots , k^{n-1} \}$ implies the special functions defined above are unique. 

\begin{thm} \label{thm:mainevent}
Let $p(x)$ be a monic $n^{th}$-order polynomial over $\Acal$ and let $D = d/d\zeta$ denote the derivative operator on $\Acal$. The special functions of the exponential for the characteristic extension algebra form a fundamental solution set to $p(D)[\eta]=0$ on $\Acal$.
\end{thm}

\noindent
{\bf Proof:} let $p(x) \in \Acal[x]$ be a monic $n$-th order polynomial; $p(x) = x^n+a_{n-1}x^{n-1}+ \cdots + a_1x+a_0$. Let $\Acal'$ be the characteristic extension algebra for $p(D)[\eta]=0$ with standard $\Acal$-module basis $1,k, \dots , k^{n-1}$.  We repeatedly apply Corollary \ref{coro:exponAprime} in what follows:
\begin{align}
p(D)[e^{k\zeta}] 
&= D^n[e^{k\zeta}]+a_{n-1}D^{n-1}[e^{k\zeta}]+ \cdots + a_1D[e^{k\zeta}]+a_0[e^{k\zeta}] \\ \notag
&= k^ne^{k\zeta}+a_{n-1}k^{n-1}e^{k\zeta}+ \cdots + a_1ke^{k\zeta}+a_0e^{k\zeta} \\ \notag
&= p(k)e^{k\zeta}.
\end{align}
But, $p(k)=0$ by construction of $\Acal'$ and it follows $p(D)[e^{k\zeta}]=0$. However,  we also have
\begin{equation}
p(D)[e^{k \zeta}] = p(D)[f_1( \zeta)]+p(D)[f_2(\zeta)]k+ \cdots + k^{n-1}p(D)[f_n(\zeta)]
\end{equation}
thus $p(D)[f_i]=0$ for $i=1,2,\dots , n$. It remains to show $\{f_1,\dots , f_n \}$ form a linearly independent set of functions. Let $h(\zeta) = e^{k\zeta}$ and note $h^{(j)}(\zeta) = k^je^{k\zeta}$ and
\begin{equation}
h^{(j)}(\zeta) = 
f_1^{(j)}(\zeta)
+f_2^{(j)}(\zeta)k
+ \cdots + f_n^{(j)}(\zeta)k^{n-1}
\end{equation}
thus $h^{(j)}(0) = k^j = f_1^{(j)}(0)
+ \cdots + f_{j+1}^{(j)}(0)k^j
+ \cdots + f_n^{(j)}(0)k^{n-1}$ for $j=0,1,\dots, n-1$. It follows that $f^{(j)}_i(0) = \delta_{i,j+1}$ for $1 \leq i, j+1\leq n$.  Consequently, we find the Wronskian of $f_1, \dots , f_n$ at $\zeta=0$ amounts to calculating the determinant of the identity matrix and hence $\{ f_1, \dots, f_n \}$ is linearly independent on $\Acal$. $\Box$ 

\begin{exa}
Consider $y''-2y'+y=0$ over $\Acal=\RN$. Identify $p(x)=(x-1)^2$ and $\Acal' = \{ a+bk \ | \ (k-1)^2=0 \}$. Thus calculate special functions from the exponential,
\begin{equation}
e^{kt} = e^{t+(k-1)t} = e^te^{(k-1)t} = e^t[1+(k-1)t] = e^t(1-t)+k(te^t)
\end{equation}
Thus $e^t(1-t), te^t$ forms a fundamental solution set for $y''-2y'+y=0$.
\end{exa}

\begin{exa}
Consider $\frac{d^n\eta}{d\zeta^n}=0$ in $\Acal$. Identify $p(x)=x^n$ gives characteristic extension algebra $\Acal'$ with basis $1,k, \dots , k^{n-1}$ where $k^n=0$. Hence calculate:
\begin{equation} 
e^{\zeta k} = 1+\zeta k+ \frac{1}{2}\zeta^2 k^2+ \cdots +  \frac{1}{(n-1)!}\zeta^{n-1}k^{n-1}. 
\end{equation}
Thus identify find special functions $1, \zeta, \zeta^2/2, \dots, \zeta^{n-1}/(n-1)!$. 
\end{exa}

\noindent
It is not always practical to calculate the exponential directly. If an isomorphism from $\Acal'$ to a well-known direct product can be found then we gain further insight. 

\begin{thm} \label{thm:isomorphismexpfla}
If $\Psi: \Acal \rightarrow \Bcal$ is an isomorphism of finite dimensional associative, commutative algebras over $\RN$ then
$ \text{exp}(\zeta) = \Psi^{-1} \left[ \text{exp}( \Psi (\zeta)) \right]$. \end{thm}

\noindent
We omit proof of the above and instead illustrate its application.

\begin{exa}
Recall $\Hcal \approxeq \RN \times \RN$ via the isomorphism $\Psi(x+jy) = (x+y,x-y)$ which has $\Psi^{-1}(a,b) = \frac{1}{2}(a+b)+\frac{j}{2}(a-b)$. Calculate:
\begin{align} 
exp( x+jy) &= \Psi^{-1}( \text{exp}( x+y,x-y)) \\ \notag
&= \Psi^{-1}(e^{x+y},e^{x-y}) \\ \notag
&= \frac{1}{2}\left( e^{x+y}+e^{x-y}\right) + \frac{j}{2}\left( e^{x+y}-e^{x-y}\right) \\ \notag
&= e^x\cosh(y)+je^x\sinh(y).
\end{align}
\end{exa}

\begin{exa}
Let us return to our original motivation for this technique. To solve $y''+y=0$ over $\RN$ we naturally consider $\Acal' = \{ a+bk \ | \ k^2=-1 \}$. Of course, $\Psi(a+bk)=a+bi$ provides the obvious isomorphism with $\CN$ and we find:
\begin{equation}
 e^{kt} = \Psi^{-1}(e^{it}) = \Psi^{-1}( \cos t+ i \sin t) = \cos t+k \sin t.
\end{equation}
Hence $y=c_1\cos t + c_2 \sin t$ is the general solution.
\end{exa}

\noindent
Notice we can just as well solve problems where the characteristic equation has real roots.

\begin{exa}
Consider $y''-3y'+2y=0$ over $\Acal= \RN$. Identify $p(x)=x^2-3x+2$ hence the characteristic extension $\Acal'$ has typical element $a+bk$ with $k^2-3k+2=0$. This algebra is isomorphic to the direct product algebra $\RN \times \RN$. In particular,
\begin{equation}
\Psi(a+bk) = (a+b, a+2b) \qquad \text{with} \qquad \Psi^{-1}(u,v) = 2u-v+(v-u)k.
\end{equation}
Therefore,  by Theorem \ref{thm:isomorphismexpfla} we calculate 
\begin{equation}
e^{kt} = \Psi^{-1}( \text{exp}(t,2t)) 
= \Psi^{-1}( (e^t,e^{2t})) 
=  2e^t-e^{2t}+(e^{2t}-e^t)k.  
\end{equation}
Thus, by Theorem \ref{thm:mainevent}, we conclude $2e^t-e^{2t},e^{2t}-e^t$ is a fundamental solution set for $y''-3y'+2y=0$.
\end{exa}

\noindent
Next we provide a less trivial application, we return solve the differential equation from Example \ref{exa:hyperbolicwithzddiffroot} by our extension technique.
 
\begin{exa}
Consider $(D-1)(D-j)[\eta]=0$ where $D = d/dz$ for hyperbolic variable $z=x+jy$. Observe the characteristic extension algebra has the form:
\begin{equation}
\Hcal' = \{ a+bk \ | \ a,b \in \Hcal, \ k^2-(1+j)k+j=0  \}.
\end{equation}
Observe $\Hcal'$ is isomorphic to the direct product of $\Gamma = \{ x+ \varepsilon y \ | \ x,y \in \RN, \varepsilon^2=0 \}$ and $\Hcal$ by the isomorphism
\begin{equation}
\Psi(a_0+ja_1+(b_0+jb_1)k) = ( a_0+a_1+(1+\varepsilon)(b_0+b_1), a_0-a_1+j(b_0-b_1))
\end{equation}
which has inverse
\begin{equation}
\Psi^{-1}(X+\varepsilon Y, U+jV) =
\frac{1+j}{2}(X-Y)+\frac{1-j}{2}U+ \left[ \frac{1+j}{2}Y+\frac{1-j}{2}V \right]k.
\end{equation}
Consider, $kz = k(x+jy)$ hence identify  $a_0=a_1=0$ and $b_0=x, b_1=y$ thus
\begin{align} 
e^{kz} &= \Psi^{-1}[ e^{( (1+\varepsilon)(x+y), j(x-y))}] \\ \notag
&=\Psi^{-1}[( e^{(1+\varepsilon)(x+y)}, e^{j(x-y)})] \\ \notag
&=\Psi^{-1}[( e^{x+y}+\varepsilon(x+y) e^{x+y}, \cosh(x-y)+j\sinh(x-y))] \\ \notag
&= \tfrac{1+j}{2}e^{x+y}(1-x-y) +\tfrac{1-j}{2}\cosh(x-y) + \left[ \tfrac{1+j}{2}(x+y) e^{x+y}+\tfrac{1-j}{2}\sinh(x-y) \right]k.
\end{align}
Thus we find fundamental solution set:
\begin{equation}
\left\{ \tfrac{1+j}{2}e^{x+y}(1-x-y) +\tfrac{1-j}{2}\cosh(x-y), \tfrac{1+j}{2}(x+y) e^{x+y}+\tfrac{1-j}{2}\sinh(x-y) \right\}.
\end{equation}
The reader may compare this to the solution derived in Example \ref{exa:hyperbolicwithzddiffroot}.
\end{exa}

\section{degenerate constant coefficient $\Acal$-ODEs} \label{sec:degenerate}

\noindent
If $\alpha \in \Acalzd$ is nonzero and $D = d/ d\zeta$ then a {\bf degenerate} $n$-th order $\Acal$-ODE has the form:
\begin{equation}
(\alpha D^n+q(D))[\eta] = 0
\end{equation}
where $q(D) \in \Acal[D]$ is at most $n-1$ degree. Degenerate $\Acal$-ODEs misbehave in many ways. We simply make a few initial observations here.

\begin{exa} \label{exa:zerdivI}
Consider the $\mathbf{\Gamma}$-ODE $\epsilon w' = 0$. Then all functions of the form $\epsilon f(z)$, where $f(z)$ is an arbitrary $\mathbf{\Gamma}$-differentiable function, are solutions.
\end{exa}

\begin{exa} \label{exa:zerdivII}
Consider the degenerate $\Hcal$-ODE $(j-1)w'=0$. Clearly $w = (j+1)f(z)+c$ solves the differential equation where $f(z)$ is any hyperbolic differentiable function and $c$ is a hyperbolic constant. The set of solutions has infinite rank as an $\Hcal$-module. Notice any initial value problem has at least one solution, but, certainly not uniquely so. For example, we may fit the initial condition $w(0)=w_0$ with the constant solution $w_1(z) = w_0$ or with any of the nonconstant solutions $w_2(z) = (j+1)z^n+w_0$ where $n \in \NN$. \end{exa}

\noindent 
McCoy's Theorem \cite{McCoy} indicates a zero divisor in a polynomial ring is annihilated by some zero divisor in the underlying ring\footnote{a nilfactorable algebra is an algebra where each zero divisor in the polynomial ring appears as a multiple of a zero divisor in the underlying ring. See \cite{bedell1} for further discussion, not all our examples are nilfactorable. }. This algebra explains the structure seen in Examples \ref{exa:zerdivI} and \ref{exa:zerdivII} and leads us to the following result:

\begin{thm}
Consider the constant coefficient $\Acal$-ODE $L[w] = 0$. This $\Acal$-ODE has solutions of the form $c f(z)$ where $c$ is a fixed constant in $\Acal$ and $f(z)$ is an arbitrary function in $C(\Acal)$ if and only if $L \in \zd(\Acal[D])$. 
\end{thm}

\noindent
{\bf Proof:} Suppose $L[cf(z)] = cL[f(z)] = 0$ for all $f(z) \in C(\Acal)$. Then $cL$ must be the zero operator. Hence, $c$ and $L$ must both be zero divisors in $\Acal[D]$. Conversely, suppose that $L \in \zd(\Acal[D])$, then by McCoy's Theorem, there exists $c \in \zd(\Acal)$ such that $cL = 0$. Hence, if $f(z) \in C(\Acal)$ is an arbitrary function, then we have $L[cf(z)] = cL[f(z)] = 0$. $\Box$ \\

\noindent There are also more subtle ways an $\Acal$-ODE might be degenerate, for example: 

\begin{exa}
Consider the direct product of $\RN$ with itself; $\Acal = \RN \times \RN$. The standard basis $e_1=(1,0), e_2=(0,1)$ serve as zero divisors in this algebra; $e_1e_2=0$. On the other hand, $e_1+e_2 = (1,1)$ is the identity of $\Acal$.  Observe that
\begin{equation}
(e_2D^2+e_1D+e_1+e_2)[\eta]=0 
\end{equation}
is a degenerate $\Acal$-ODE. We can solve it by direct calculation. Since $D=d/d\zeta = ( \partial_x, \partial_y)$ and the $\Acal$-CR equations for $\eta \in \FunA$ provide $\eta = (\eta_1, \eta_2)$ has $\partial_y \eta_1 = 0$ and $\partial_x \eta_2 = 0$. Consequently, we find:
\begin{equation}
 (\partial_x+1, \partial_y^2+1)(\eta_1, \eta_2) = (0,0) \ \ \Rightarrow \ \
\frac{d\eta_1}{dx}+\eta_1 = 0, \ \  \frac{d^2\eta_2}{dy^2}+\eta_2 = 0.
\end{equation}
Therefore, the general solution is simply:
$
\eta(x,y) = (c_1e^{-x}, c_2 \cos y+c_3 \sin y)
$.  Yet, if we impose the initial value condition $\eta(0,0) = (0,0), f'(0,0) = (1,1)$, we see that this initial value problem has no solution.
\end{exa}

\noindent 
In fact, this behavior is typical. If a degenerate $\Acal$-ODE $L[\eta]=0$ has $L \notin \zd( \Acal[D])$ then there exist initial value problems for $L[\eta]=0$ which have no solution. This contrasts Example \ref{exa:zerdivII} where every initial value problem had a (non-unique) solution.

\begin{thm}
Suppose $\alpha \in \Acalzd$ is nonzero and $p(D) = \alpha D^n+q(D)$ where $q(D) \in \Acal[D]$ is at most $n-1$ degree and $q(D)$ is not a multiple of $\alpha$. Then for any $\zeta_o \in \Acal$, there exists $c_0, c_1, \dots , c_{n-1} \in \Acal$ such that no solution $\eta$ to $p(D)[\eta]=0$ has
$ \eta(\zeta_o) = c_0, \eta'(\zeta_o)=c_1, \dots, \eta^{(n-1)}(\zeta_o)= c_n$. 
\end{thm}

\noindent
{\bf Proof:} Suppose $\alpha \in \Acalzd$ is nonzero and $p(D) = \alpha D^n+q(D)$ where $q(D) = q_{n-1}D^{n-1}+ \cdots + q_1D+q_0$ where $q(D)$ is not a multiple of $\alpha$. Hence there exists $q_j$ with $0 \leq j \leq n-1$ for which $q_j$ is not a multiple of $\alpha$. Let $\zeta_o \in \Acal$ and suppose we wish to fit the data $\eta^{(i)}(\zeta_o) = \delta_{ij}$ meaning $\eta^{(i)}(\zeta_o) = 0$ for $i=0,1,\dots , n-1$ where $i \neq j$ and $\eta^{(j)}(\zeta_o)=1$. Observe such a solution has $q(D)[\eta](0) = q_j$ and hence $p(D)[\eta](0) = 0$ implies $\alpha \eta^{(n)}(0) - q_j = 0$. Thus $q_j = \alpha \eta^{(n)}(0)$ and so $q_j$ is a multiple of $\alpha$. Hence, using proof by contradiction, we conclude there is no solution to $p(D)[\eta]=0$ with $\eta^{(i)}(\zeta_o) = \delta_{ij}$. $\Box$ 

\section{Cauchy Euler problems} \label{sec:cauchyeuler}

A particularly simple class of non-constant coefficient $\mathcal{A}$-ODEs are the so-called Cauchy Euler problems, which are $\mathcal{A}$-ODEs $L[w] = f$ where $L$ is a differential operator of the form $ a_n z^n D^n + \dots + a_1 z D + a_0 $. To avoid degeneracies as in Section \ref{sec:degenerate} we assume $a_n \in \Acalx$ and thus set $a_n=1$ without loss of generality in what follows.  \\

\noindent
Cauchy Euler problems require we pay more attention to domain than in the constant coefficient case. If we suppose $z \in \mathrm{Ld}(\mathcal{A})$ then the function $z^\alpha = e^{\alpha \log(z)}$ is well-defined for all $\alpha \in \Acal$. Consult \cite{bedell2} for details on the construction of the logarithm including the structure of its domain $\mathrm{Ld}(\mathcal{A})$. Furthermore, a necessary and sufficient condition that $w = z^\alpha$ solves $L[w] = 0$ on $\mathrm{Ld}(\mathcal{A})$ is given by the {\bf characteristic equation}
\begin{equation}
  \alpha (\alpha - 1) \dots (\alpha - n + 1) + \dots + a_2\alpha(\alpha -1)+a_1 \alpha +a_0= 0. 
\end{equation}
It is also important to note that even though we assumed that $z \in \Ld(\Acal)$ in order to formulate the problem in general, if the roots of the characteristic polynomial are integers, we may be able to obtain better results, for example:

\begin{exa}
Consider the Cauchy Euler problem $z^2 w'' + z w' - w = 0$ in $\Hcal$. Supposing solutions are of the form $w = z^\alpha$ we obtain characteristic equation $\alpha^2 - 1 = (\alpha - 1)(\alpha + 1) = 0$. Hence, $z^1$ and $z^{-1} = \frac{1}{z}$ both give solutions to the Cauchy Euler problem. However, $z$ gives a solution on all of $\Acal$, and $z^{-1}$ gives a solution on all of $\Acal^\times$.
\end{exa}

\noindent
Similar cases occur for the usual Cauchy Euler problem on $\RN$ where polynomial solutions may exist on $\RN$ whereas other solutions merely exist on $(0,\infty)$. Noting that positive integer power functions $z \mapsto z^n$ are entire, while $z \mapsto \frac{1}{z^n}$ are $\Acal$-differentiable on $\Acalx$ we find:


\begin{prop}
Let $z^n w^{(n)} + \dots + a_1 z w' + a_0 w = 0$ be a Cauchy Euler problem over the algebra $\Acal$ and let $p(\alpha)$ be its characteristic equation. Then if $\alpha$ is a root of $p(\alpha)$:
\begin{quote}
 \begin{enumerate}[{\bf (i.)}]
\item If $\alpha \in \mathbb{N}$, $z^\alpha$ is a solution to the Cauchy Euler problem on all of $\Acal$.
\item If $\alpha \in \mathbb{Z}^-$, $z^\alpha$ is a solution to the Cauchy Euler problem on $\Acal^\times$.
\item Otherwise, $z^\alpha$ is a solution to the Cauchy Euler problem on $\Ld(\Acal)$.
\end{enumerate}
\end{quote}
\end{prop}

\noindent
The result above is not general. For example, 
\begin{equation} \label{eqn:noCEroots}
 z^2 w'' + (j+1) z w' + j w = 0  
\end{equation} 
in $\mathcal{H}$, has characteristic equation $\alpha^2 + j \alpha + j = 0$ for which no solutions over $\Hcal$ exist. Thus solutions  of the form $z^{\alpha}$ cannot be found for Equation \ref{eqn:noCEroots}.  Repeated roots in the characteristic equation also pose additional problems. However, we can solve both of these issues by applying a substitution that transforms the Cauchy Euler problem into an constant coefficient $\Acal$-ODE. In particular, consider 
\begin{equation} \label{eqn:cauchyeulerprob}
 z^n w^{(n)} + \dots + a_1 z w' + a_0 w = 0
\end{equation} 
in the algebra $\Acal$ and make the substitution of variables $z = e^\zeta$ given that $z \in \Ld(\Acal)$. Thus, if $\log : \Ld(\Acal) \rightarrow B$ is a branch of the logarithm, then we also have $\zeta = \log(z)$. We also define $\eta(\zeta) = w(e^\zeta)$ which is to say $w(z) = \eta( \log(z))$. From this we have: 
\begin{equation} \label{eqn:chooseI}
 \frac{dw}{dz} = \zeta'(z)\eta'(\zeta)  = \frac{1}{z} \eta'(\zeta)  
 \end{equation}
Continue differentiating to derive:\begin{equation} \label{eqn:chooseII}
 \frac{d^2w}{dz^2} = \frac{d}{dz} \left[
 \frac{1}{z}\eta'(\zeta) \right]  = \frac{-1}{z^2}\eta'(\zeta)+ \frac{1}{z}\frac{d}{dz}\left[ \eta'(\zeta) \right] = \frac{1}{z^2}\left[\eta''(\zeta)-\eta'(\zeta) \right]
 \end{equation}
and,
\begin{equation} \label{eqn:chooseIII}
 \frac{d^3w}{dz^3} = \frac{d}{dz} \left[
\frac{1}{z^2}\left[\eta''(\zeta)-\eta'(\zeta) \right] \right]  = \frac{1}{z^3}\left[\eta'''(\zeta)-3\eta''(\zeta)+2\eta'(\zeta)  \right]
 \end{equation}
We note $z^k \frac{d^kw}{dz^k} = \frac{d^k\eta}{d\zeta^k}+ q(d/d\zeta)[\eta]$ where $q(D) = q_{k-1}D^{k-1}+ \cdots + q_1D+q_0$ for appropriate coefficients $q_{k-1}, \dots , q_1,q_0$. Therefore, under the substitution $z = e^\zeta$ for $z \in \Ld(\Acal)$ and $w(z) = \eta( \log(z))$ we find the Cauchy Euler problem of Equation \ref{eqn:cauchyeulerprob} transforms to:
\begin{equation}
\frac{d^k\eta}{d\zeta^k}+b_{k-1}\frac{d^{k-1}\eta}{d\zeta^{k-1}}+ \cdots +b_1\frac{d\eta}{d\zeta}+b_0\eta = 0
\end{equation}
for appropriate constants $b_{k-1}, \dots , b_1, b_0 \in \Acal$. Equations \ref{eqn:chooseI}, \ref{eqn:chooseII} and \ref{eqn:chooseIII} indicate how to transform a Cauchy Euler problem of order $n \leq 3$. Therefore, we can find a fundamental solution set for the transformed problem and hence a fundamental solution set for the Cauchy Euler problem on $\Ld(\Acal)$. 

\begin{exa}
Observe $z^2w''+3zw'+w=0$ with $z \in \Ld(\Acal)$ transforms to 
\begin{equation}
\eta''-\eta'+3(\eta')+\eta = \eta''+2\eta'+\eta=0 
\end{equation} 
by Equations \ref{eqn:chooseI} and \ref{eqn:chooseII}. Using methods developed earlier in this article, we find general solution $\eta(\zeta) = c_1e^{-\zeta}+c_2\zeta e^{-\zeta}$ thus as $w(z) = \eta( \log(z))$ we deduce: using $z^{\alpha} = e^{\alpha \log(z)}$ for $z \in \Ld( \Acal)$,
\begin{equation}
w(z) = \frac{c_1}{z}+ \frac{c_2\log(z)}{z}.
\end{equation} 
\end{exa}

\noindent
In the example above we intentionally left $\Acal$ as arbitrary. The structure of the logarithm is explicitly quite different for various choices of $\Acal$. 

\begin{exa}
Consider the Cauchy Euler problem over $\Hcal = \RN \oplus j \RN$ where $j^2=1$:
\begin{equation} \label{eqn:CEproblemII}
z^2w''+zw'-\frac{1}{4}w = 0
\end{equation}
If $z=x+jy \in \Ld( \Hcal)$ then there exists $\zeta \in \Hcal$ for which $z = e^{\zeta}$ and we denote $\log(z) = \zeta$ for such $z$. In particular, $x+jy \in \Ld( \Hcal)$ requires $x>0$ and $-x< y< x$ and we calculate:
\begin{equation}
\log(x+jy) = ln \sqrt{x^2-y^2}+j \tanh^{-1}(y/x).
\end{equation}
Transforming Equation \ref{eqn:CEproblemII} by Equations \ref{eqn:chooseI} and \ref{eqn:chooseII} yields
\begin{equation}
\eta''-\frac{1}{4}\eta = 0 \ \ \Rightarrow \ \ \eta(\zeta) = c_1e^{-\zeta/2}+c_2e^{-\zeta/2} \ \ \Rightarrow \ \ w(z) = c_1e^{-\frac{1}{2}\log(z)}+c_2e^{\frac{1}{2}\log(z)}
\end{equation}
or, using $z^{\alpha} = e^{\alpha \log(z)}$ we find general solution $w(z) = c_1z^{-1/2}+c_2z^{1/2}$ for $z \in \Ld( \Hcal)$. Explicitly, after a little calculation,
\begin{equation} \label{eqn:hyperbolicroot}
(x+jy)^{1/2} = \frac{1}{\sqrt{2}}\sqrt{ x+ \sqrt{x^2-y^2}}+\frac{j}{\sqrt{2}}\sqrt{x-\sqrt{x^2-y^2}}
\end{equation}
For example, $(5+4j)^{1/2} = 2+j$.
\end{exa}
 
\section{acknowledgements}
The authors are thankful to R.O. Fulp for helpful comments on a rough draft of this article.

\end{document}